\documentclass[12pt, letterpaper]{elsarticle}
\usepackage{amsmath}
\usepackage{amssymb}
\usepackage{graphicx}
\usepackage{multirow}
\usepackage{latexsym}
\usepackage{epic}
\usepackage{url}
\usepackage{soul}
\usepackage{afterpage}
\usepackage{comment}
\usepackage{bbm}
\usepackage{setspace}
\usepackage[norelsize,ruled, lined,linesnumbered]{algorithm2e}
\usepackage{enumitem}
\usepackage{empheq,mathtools}
\usepackage[symbol]{footmisc}
\usepackage[mathscr]{euscript}
\usepackage{tikz}
\usepackage{tikz,fullpage}
\usepackage{multirow}
\usepackage{pgf}
\usetikzlibrary{arrows,automata}
\usepackage{tkz-berge}
\usepackage{amsthm}
\usepackage{subcaption}
\usepackage[toc,page]{appendix}

\usepackage[left=0.9in,top=0.9in,right=0.9in,bottom=0.9in,nohead]{geometry}

\makeatletter
\newtheorem*{rep@theorem}{\rep@title}
\newcommand{\newreptheorem}[2]{%
\newenvironment{rep#1}[1]{%
 \def\rep@title{#2 \ref{##1}}%
 \begin{rep@theorem}}%
 {\end{rep@theorem}}}
\makeatother

\newreptheorem{example}{Example}
\newtheorem{theorem}{Theorem}

\newtheorem{proposition}{Proposition}

\newtheorem{example}{Example}

\usepackage{xcolor}
\usepackage[colorinlistoftodos,prependcaption,textsize=tiny]{todonotes}
\usepackage{xargs}

\usepackage[symbol]{footmisc}

\begin{document}
\onehalfspace
\begin{frontmatter}
\title{On the multi-stage shortest path problem under distributional uncertainty}

\author[label1]{Sergey S.~Ketkov\footnote[2]{Corresponding author. Email: sketkov@hse.ru; phone: +7-910-382-2732.}}

\address[label1]{Laboratory of Algorithms and Technologies for Networks Analysis, HSE University,\\ Rodionova st., 136, Nizhny Novgorod, 603093, Russia}

\begin{abstract}
In this paper we consider an ambiguity-averse multi-stage network game between a \textit{user} and an \textit{attacker}. The arc costs are assumed to be random variables that satisfy prescribed first-order moment constraints for some subsets of arcs and individual probability constraints for some particular arcs. The user aims at minimizing its cumulative expected loss by traversing between two fixed nodes in the network, while the attacker's objective is to maximize the user's expected loss by selecting a distribution of arc costs from the family of admissible distributions. In contrast to most of the related studies, both the user and the attacker can dynamically adjust their decisions at particular nodes of the user's path. By observing the user's decisions, the attacker may reveal some additional distributional information associated with the arcs emanated from the current user's position. It is shown that the resulting multi-stage distributionally robust shortest path problem (DRSPP) admits a linear mixed-integer programming reformulation (MIP). 
In particular, we distinguish between acyclic and general graphs by introducing different forms of non-anticipativity constraints. Finally, we perform a numerical study, where the quality of adaptive decisions and computational tractability of the proposed MIP reformulation are explored with respect to several classes of synthetic network~instances.
\end{abstract}

\begin{keyword}
shortest path problem; distributionally robust optimization; polyhedral uncertainty; piecewise constant decision rules; mixed-integer programming
\end{keyword}

\end{frontmatter}
\onehalfspace
\section{Introduction} \label{sec: intro}
\textit{Distributionally robust optimization} (DRO) is a methodology for addressing uncertainty in optimization problems, where the probability distribution of uncertain parameters is only known to reside within an \textit{ambiguity set} (or a \textit{family}) of admissible distributions; see, e.g., the studies in \cite{Bertsimas2018, Delage2010, Esfahani2018, Goh2010, Wiesemann2014}. A standard one-stage DRO problem can be usually formulated as a bi-level problem, where the expected value of the loss function (or some measure of risk, if the decision-maker is risk-averse) is minimized under the worst-case possible distribution of uncertain parameters. 

Typically, in one-stage DRO problems the decisions are implemented \textit{here-and-now}, before the realization of uncertainty. However, 
the distributionally robust optimization approach may also address dynamic (or, equivalently, multi-stage) optimization problems, where decisions adapt to the uncertain outcomes as they unfold in stages; see, e.g., \cite{Bertsimas2018, Goh2010, Hanasusanto2016, Jiang2018}. In other words, in multi-stage DRO problems decision variables at the current stage can be thought as some functions of the uncertain problem parameters observed up to this stage.

In this study we consider a multi-stage version of the \textit{distributionally robust shortest path problem} (DRSPP), where the vector of arc costs/travel times in a given network is subject to uncertainty. From the game theoretical perspective, our problem can be viewed as a dynamic zero-sum game between two decision-makers, which are referred to as a \textit{user} and an \textit{attacker}. The user attempts to minimize its expected loss by traversing between two fixed nodes in a given network. On the other hand, the attacker aims at maximizing the user's expected loss by selecting a distribution of arc costs within a given family of probability distributions. The outlined game is dynamic in the sense that both the user and the attacker are able to adjust their decisions at particular nodes of the user's path.

In fact, our multi-stage formulation is motivated and built upon the related one-stage formulation of DRSPP in \cite{Ketkov2021}. Specifically, we preserve the same form of the ambiguity set and the objective function but introduce some \textit{auxiliary distributional constraints} that can be verified by the user dynamically while traversing through the network. 
 
\subsection{Related literature} \label{subsec: related literature}

 A rather large number of studies consider one-stage versions of DRSPP assuming different forms of ambiguity sets and objective criteria; see, e.g., \cite{Cheng2016, Gavriel2012, Wang2020, Zhang2017b}. 
Under some assumptions on the geometry of ambiguity sets and functional properties of the objective function, one-stage problems can be reduced to single-level linear or non-linear mixed-integer programming (MIP) problems. In particular, these single-level reformulations typically rely on strong duality results for moments problems due to Isii \cite{Isii1962} and Shapiro \cite{Shapiro2001}. 
We also refer to \cite{Khanjani2021} for another version of DRSPP, in which the uncertain problem parameters are related not to the arc costs/travel times in a given network but to some resource constraints.

It can also be argued that in a number of practical applications the arc costs/travel times in a given network can only be observed using a finite training data set; see, e.g.,~\cite{Cao2014,Chassein2019,Wang2020}. In this case one may attempt to construct an ambiguity set, which contains the nominal distribution of the cost vector with a required confidence level. For example, some data-driven approaches to constructing ambiguity sets are based on estimating the mean and the covariance matrix of uncertain problem parameters \cite{Delage2010} or using some distance metrics from the empirical distribution of the data \cite{Bayraksan2015, Esfahani2018}.

\looseness-1 In contrast to one-stage problems, multi-stage distributionally robust formulations of the shortest path problem are considered by relatively few authors. The major limitation of this problem setting is that \textit{recourse decisions}, i.e., decisions affected by uncertainty, are discrete and, therefore, the application of standard linear or piecewise linear decision rules \cite{Bertsimas2018, Goh2010} is rather limited. In the following, we discuss multi-stage problem formulations that may account the shortest path problem as a special case. 

First, we refer to the study by~Hanasusanto~et~al.~\cite{Hanasusanto2016}, where a two-stage DRO problem with binary recourse decisions is approximated by using a restricted number of preselected second-stage decisions. The authors prove that the proposed approximation is exact, if the number of preselected decisions is sufficiently large. Despite the fact that the study in \cite{Hanasusanto2016} may address rather general types of ambiguity sets and objective criteria, it is not quite clear how to extend these theoretical results to a multi-stage problem setting.

Another line of research is focused on exact dynamic programming based algorithms assuming that the uncertain parameters are \textit{stage-wise independent}. To the best of our knowledge, Yu and Shen~\cite{Yu2020} are among the first, who consider multi-stage distributionally robust MIPs. The authors investigate the concept of decision-dependent ambiguity sets, where parameters of the first- and the second-order moment constraints at some stage $t \in \mathbb{Z}_+$ are linearly dependent on the user's decision at the previous stage, $t-1$. The authors reformulate the stage-wise problems as linear MIP or semidefinite programming problems. Then different versions of stochastic dual dynamic integer programming (SDDiP) algorithm \cite{Zou2019} are applied for solving the resulting multi-stage problems exactly or deriving objective bounds. The latter algorithm is shown to be finitely convergent and rather effective for moderate-sized problem instances.

In conclusion, we discuss a multi-stage shortest path interdiction model proposed by Sefair~et~al.~\cite{Sefair2016} and a multi-stage robust MIP problem formulation of Bertsimas and Dunning~\cite{Bertsimas2016}. 
Despite the fact that these models are deterministic, they provide some interesting insights for the multi-stage DRSPP considered in the current study. 

Sefair~et~al.~\cite{Sefair2016} consider a zero-sum multi-stage network game between a user and an attacker, where the attacker attempts to maximize the user's loss by blocking a subset of arcs (subject to some budgetary constraint) any time the user reaches a node in the network. Furthermore, the user may dynamically alter its path by observing the arcs blocked by the attacker. The authors in \cite{Sefair2016} demonstrate that the attacker's problem can be solved in polynomial time for acyclic graphs and admits an exact exponential-state dynamic-programming algorithm in the general case. 

On the other hand, Bertsimas and Dunning~\cite{Bertsimas2016} consider a class of linear MIP problems, where the uncertain problem parameters are only known to reside within a polyhedral uncertainty set. The uncertainty set is partitioned into a number of disjoint polyhedrons and the recourse decisions are assumed to be piecewise constant functions on the generated subsets. This modeling approach enables to obtain approximations of multi-stage robust MIP problems, where the quality of approximation is adjusted by the choice of an appropriate partition scheme. 

\subsection{Our approach and contributions} \label{subsec: approach and contributions}
In this study we design a multi-stage version of DRSPP, which is based on the related one-stage model in \cite{Ketkov2021} but allows to capture a successive revelation of distributional information to the user. In this regard, we first discuss our construction of the ambiguity set and the objective function for the one-stage problem. Next, we formulate some research questions that need to be addressed when the problem unfolds in stages. 

\textbf{One-stage problem.}
We consider a weighted connected directed graph $G := (N, A, \mathbf{c})$, where $N$ and $A$ denote its sets of nodes and directed arcs, respectively, and $\mathbf{c} = \{c_a, a \in A\}$ is a nonnegative cost vector. We also denote by $s \in N$ and $f \in N$ the source and the destination nodes in $G$, respectively. Then a standard \textit{risk-neutral one-stage} DRSPP can be formulated as: 
\begin{equation} \label{intro: single-stage DRO}
\min_{\mathbf{y} \in Y} \max_{\mathbb{Q} \in \mathcal{Q}} \mathbb{E}_{\mathbb{Q}} \{\ell(\mathbf{c}, \mathbf{y})\},
\end{equation}
 where $Y \subseteq \{0,1\}^{|A|}$ is a set of all incidence vectors corresponding to feasible $s-f$ paths in $G$, $\mathbf{y} \in Y$ is an incidence vector of a path selected by the user 
 and $\ell(\mathbf{c}, \mathbf{y})$ is a given loss function. Furthermore, for any fixed $\mathbf{y} \in Y$ the attacker attempts to maximize the user's expected loss by selecting a distribution $\mathbb{Q}$ of the cost vector $\mathbf{c}$ from some family of distributions $\mathcal{Q}$.

In the related one-stage model in \cite{Ketkov2021} 
we assume that the user has the following \textit{partial information} about the distribution $\mathbb{Q}$.
First, the costs of particular arcs $a \in A$ are subject to \textit{individual probability constraints} of the form:
\begin{equation} \label{probability constraints}
\mathbb{Q}\{c_a \in [l^{(j)}_a, u^{(j)}_a]\} \in [\underline{q}^{(j)}_a, \overline{q}^{(j)}_a] \quad \forall j \in \mathcal{D}_a, \mbox{ } a \in A,
\end{equation}
where $[l^{(j)}_a, u^{(j)}_a] \subseteq [l_a, u_a]$, $j \in \mathcal{D}_a = \{1, \ldots, d_a\}$, is a finite set of subintervals; $ \underline{q}^{(j)}_a$ and $\overline{q}^{(j)}_a$ are contained in $[0, 1]$ and bound the probability that the random cost $c_a$ belongs to the $i$-th subinterval. In particular, we assume that for each $a \in A$ and $j = 1$ the constraints
$$\mathbb{Q}\{c_a \in [l^{(j)}_a, u^{(j)}_a]\} \in [\underline{q}^{(j)}_a, \overline{q}^{(j)}_a]$$
are \textit{support constraints} with $l^{(j)}_a = l_a$, $u^{(j)}_a = u_a$ and $\underline{q}^{(j)}_a = \overline{q}^{(j)}_a = 1$ for some $l_a, u_a \in \mathbb{R}_{>0} \cup \{0\}$. 
Secondly, for a real-valued matrix $\mathbf{B} \in \mathbb{R}^{|A|\times k}$ and a vector $\mathbf{b} \in \mathbb{R}^{k}$, $k \in \mathbb{Z}_+$, we introduce \textit{linear expectation constraints} of the form:
\begin{equation} \label{linear expectation constraints}
\mathbb{E}_{\mathbb{Q}} \{\mathbf{B} \mathbf{c}\} \leq \mathbf{b}.
\end{equation}
Hence, the ambiguity set $\mathcal{Q}$ is formed by all distributions that satisfy the distributional constraints (\ref{probability constraints}) and (\ref{linear expectation constraints}).

We argue that in contrast to the moment-based ambiguity sets \cite{Delage2010} or ambiguity sets based on a distance metric from the empirical distribution \cite{Bayraksan2015, Esfahani2018}, our distributional constraints allow incomplete knowledge of the training data set. Formally, the individual probability constraints (\ref{probability constraints}) can be constructed by using \textit{interval-censored observations} with respect to particular arcs; see, e.g., \cite{Sun2007}. From the application perspective, the interval uncertainty can be motivated by measurement errors or sensors detection limits in environmental measurements; see, e.g., \cite{Kreinovich2006}. At the same time, the linear expectation constraints (\ref{linear expectation constraints}) can be constructed by leveraging \textit{linear combinations of arc costs} with respect to some subsets of arcs. 
For example, in our problem setting the user may only have access to random observations of a total cost with respect to some routes in the network; this type of feedback is also known as ``bandit'' feedback in online learning problem settings~\cite{Bubeck2012}. We refer the reader to Section 2.2 in \cite{Ketkov2021} and Section \ref{subsec: constraints form data} in the current study for further details on the construction of distributional constraints from data. 

The key theoretical result of \cite{Ketkov2021} indicates that the one-stage DRSPP of the form (\ref{intro: single-stage DRO}) with a linear loss function $\ell(\mathbf{c}, \mathbf{y}) = \mathbf{c}^\top \mathbf{y}$ can be reformulated as a robust shortest path problem with some polyhedral uncertainty~set; see, e.g., \cite{Buchheim2018} for a review of robust combinatorial optimization problems. Furthermore, the resulting problem admits a linear mixed-integer programming (MIP) reformulation and, thus, can be solved using off-the-shelf MIP solvers. 

\textbf{Multi-stage problem.}
In the one-stage formulation the user selects a path \textit{here-and-now} before the realization of uncertainty, i.e., some unknown distribution~$\mathbb{Q} \in \mathcal{Q}$. In the \textit{multi-stage} formulation of DRSPP we, in turn, attempt to address the following research questions:
\begin{itemize}
	\item[\textbf{Q1.}] \textit{Is there a benefit for the user to alter the chosen path, if it observes some additional distributional information while traversing through the network?}
	\item[\textbf{Q2.}] \textit{How much can the user gain by leveraging such adaptive decisions?} 
	\item[\textbf{Q3.}] \textit{Can the resulting multi-stage formulation be solved at hand using off-the-shelf MIP solvers?}
\end{itemize}

We address a dynamic revelation of distributional information to the user by introducing some \textit{auxiliary distributional constraints} (in addition to the constraints in $\mathcal{Q}$ that are known to the user a priori) associated with the arcs emanated from some current user's position. We assume that the auxiliary constraints can be \textit{verified} by the user while traversing through the network. Simply speaking, the user forms a list of auxiliary constraints \textit{at the beginning of the game}. Then for each constraint in the list the attacker decides whether this constraint is satisfied or not and reveals its response, i.e., ``yes'' or ``no'', to the user \textit{as soon as} the user achieves the respective node in the network. 

\looseness-1 From the practical perspective, the auxiliary distributional constraints can be motivated, e.g., by some additional data observed by the user from Bluetooth sensors \cite{Haghani2010, Asudegi2013} or GPS-equipped floating vehicles \cite{Vo2011}. In other words, if sensors are placed at particular nodes of the given network, then the user may quantify the travel times between two successive detection stations; see, e.g., \cite{Haghani2010}. This additional information can be collected dynamically while traversing through the network and used to verify the auxiliary constraints; we refer to Example \ref{example 2} and Section \ref{subsec: constraints form data} for a more comprehensive discussion. 

Following the one-stage formulation, we suppose that the user's objective is to select a path with the least possible expected cost. The attacker, in turn, attempts to maximize the user's objective function by selecting a distribution of arc costs from the ambiguity set $\mathcal{Q}$. However, in contrast to the one-stage problem, both decision-makers are able to adjust their decisions at particular nodes of the user's~path. In other words, the user's problem unfolds in stages, where at each stage the attacker reveals (if necessary) some new distributional information to the user and the user picks a subsequent node of its path. 

 To the best of our knowledge, the concept of using auxiliary distributional constraints is not discussed in the context of multi-stage DRO problems. The idea of our approach is to partition the initial ambiguity set into a number of disjoint subsets, where each subset corresponds to some vector of attacker's responses (recall that the attacker provides a binary response to each auxiliary constraint). Then we show that the user's decision can be viewed as a \textit{piecewise constant} function on the generated ambiguity sets subject to some \textit{non-anticipativity constraints}; see, e.g., \cite{Bertsimas2016, Goh2010}. 

In contrast to the study by Sefair~et~al. \cite{Sefair2016}, we consider more general types of constraints from the attacker's perspective (instead of a unique budget constraint) and account some partial distributional information available to the user. In contrast to the study by Yu and Shen~\cite{Yu2020}, our model exploits linear expectation constraints of the form (\ref{linear expectation constraints}), which make the uncertain parameters stage-wise dependent and thereby ``destroy'' dynamic programming based algorithms. Finally, our solution approach is, in a sense, similar to the one developed by Bertsimas and Dunning \cite{Bertsimas2016}. However, from the methodological point of view, the partitions in \cite{Bertsimas2016} are not motivated by any practical applications and are used to provide a better approximation for the nominal multi-stage problem. Meanwhile, the partition of the ambiguity set in our setting is motivated by the auxiliary distributional constraints observed by the user in the considered~game.

Our contributions for the multi-stage model can be summarized as follows:
\begin{itemize}
	\item We formulate the multi-stage DRSPP and describe two classes of non-anticipativity constraints, for acyclic and general graphs, respectively (Sections \ref{subsec: modeling assumptions} and \ref{subsec: multi-stage}).  \item We show the auxiliary distributional constraints can be constructed and verified using some data from Bluetooth sensors (Section~\ref{subsec: constraints form data}). 
	\item Under some mild assumptions, the multi-stage problem is reformulated as a one potentially large linear MIP problem (Section \ref{sec: solution techniques}). These results address our research~question~\textbf{Q3}. 
	\item The obtained MIP reformulation is used in our numerical study, where its computational tractability and the quality of adaptive decisions are explored numerically (Section~\ref{sec: comp study}). Hence, we address the research questions \textbf{Q1} and \textbf{Q2}. 
\end{itemize}


\section{Base model} \label{sec: problem}
\textbf{Notation.} All vectors and matrices are labeled by bold letters. For a network $G := (N, A, \mathbf{c})$ we denote by $N$ and $A$ its sets of nodes and directed arcs, respectively, whereas $\mathbf{c}$ is a nonnegative cost vector. Let $s$ and $f$ be the source and the destination nodes, respectively. For each node $i \in N$ we refer to $RS_i$ ($FS_i$) as the set of arcs directed out of (and into) node~$i$.
The space of all probability distributions on $\mathbb{R}^{|A|}$ is denoted as $\mathcal{Q}_0(\mathbb{R}^{|A|})$. Finally, for each $a \in A$ we denote by $\mathbb{Q}_a \in \mathcal{Q}_0(\mathbb{R})$ the marginal distribution induced by some joint distribution $\mathbb{Q} \in \mathcal{Q}_0(\mathbb{R}^{|A|})$. 

\subsection{Modeling assumptions} \label{subsec: modeling assumptions}
As outlined in Section \ref{sec: intro}, we consider a directed weighted connected graph $G := (N, A, \mathbf{c})$, where the distribution $\mathbb{Q}$ of the cost vector $\mathbf{c}$ belongs to an ambiguity set $\mathcal{Q}$ given by: 
\begin{equation} \label{ambiguity set}
\begin{gathered}
\mathcal{Q}: = \Big\{\mathbb{Q} \in \mathcal{Q}_0(\mathbb{R}^{|A|}) :\; \mathbb{Q} \mbox{ satisfies the constraints (\ref{probability constraints}) and (\ref{linear expectation constraints})}
\Big\}.
\end{gathered}
\end{equation}
 Furthermore, we note that the set $Y$ of all feasible path-incidence vectors can be expressed as: 
\begin{subequations} \allowdisplaybreaks \label{path flow constraints}
	\begin{align} 
	Y = \Big\{\mathbf{y} \in \{0, 1\}^{|A|}: \; & \sum_{a \in FS_i} y_a - \sum_{a \in RS_i} y_a = \begin{cases}
	1, \mbox{ if } i = s \\
	-1 \mbox{ if } i = f \\
	0, \mbox{ otherwise } 
	\end{cases} \quad \forall i \in N \label{eq: flow conservation} \\
	& \sum_{a \in FS_i} y_a \leq 1 \quad \forall i \in N \label{eq: exlude cycles}
	\Big\},
	\end{align}
\end{subequations}
where the constraints (\ref{eq: flow conservation}) are standard flow conservation constraints and the constraints (\ref{eq: exlude cycles}) ensure that each node is visited at most once; see, e.g., \cite{Taccari2016}. We recall that negative cycles are prohibited due to the support constraints, i.e., the probability constraints~(\ref{probability constraints}) for each $a \in A$ and $j = 1$. 

Throughout the paper we make the following modeling assumptions:
\begin{itemize}
	\item[\textbf{A1.}] Each node $i \in N$ is visited by the user at most once.
	\item[\textbf{A2.}] Both the user and the attacker have complete information about the initial family of distributions $\mathcal{Q}$ and the network $G$.
	\item[\textbf{A3.}] In addition to the distributional constraints in $\mathcal{Q}$, the user forms a list $\mathcal{L}$ of auxiliary distributional constraints given by:
\begin{subequations} \label{eq: list of additional constraints} \allowdisplaybreaks
	\begin{align}
	\mathcal{L} := \bigcup_{i \in N \setminus \{f\}} \Big\{& \mathbb{Q}_a\{c_a \in [\widetilde{l}^{(j)}_a, \widetilde{u}^{(j)}_a]\} \leq \widetilde{q}^{\; (j)}_a \quad \forall j \in \widetilde{\mathcal{D}}_a, \; \forall a \in FS_i \; ; \label{cons: auxiliary probability}\\
	& \mathbb{E}_{\mathbb{Q}}\{\sum_{a \in FS_i} p_{ja} c_a\} \leq p_{j0} \quad \forall j \in \widetilde{\mathcal{K}}_i
	\Big\}, \label{cons: auxiliary expectation}
	\end{align}
\end{subequations}
where, $\widetilde{\mathcal{D}}_a = \{1, \ldots, \widetilde{d}_a\}$, $a \in A$, and $\widetilde{\mathcal{K}}_i = \{1, \ldots, \widetilde{k}_i\}$, $i \in N \setminus \{f\}$, are potentially empty sets of indexes.\end{itemize}

Assumption \textbf{A1} can be explained, e.g., by a nonstationarity in the attacker's distribution (a similar assumption is made by Xu and Mannor \cite{Xu2010} with regard to distributionally robust Markov decision processes). Put differently, if the user returns to a node multiple times, then the attacker may simply modify the underlying distribution of arc costs. In this situation the user cannot exploit the previously collected distributional information and, hence, it is not favorable for the user to visit nodes multiple times. 


 In the first part of Assumption \textbf{A2} we envision that the user has some initial information about the distribution of arc costs. As briefly outlined in Section \ref{subsec: approach and contributions}, this information can be collected beforehand by leveraging some historical data. The second part of Assumption \textbf{A2} is rather standard in the robust and distributionally robust optimization literature. 
Nevertheless, we refer to the studies in \cite{Borrero2015, Borrero2019, Ketkov2020} for multi-stage shortest path interdiction models, in which the attacker has incomplete knowledge about the structure of the underlying network and observes the existence and precise costs of particular arcs by observing the user's decisions. 

Assumption \textbf{A3} indicates that the auxiliary distributional constraints (\textit{i}) have a form similar to the initial distributional constraints in $\mathcal{Q}$ and (\textit{ii}) are associated only with the arcs emanated from some node $i \in N \setminus \{f\}$. In the following, we demonstrate that the first property, (\textit{i}), is necessary for deriving a linear MIP reformulation of the proposed multi-stage optimization problem. On the other hand, violation of the second property, (\textit{ii}), is not favorable for the attacker (as it may allow the user to observe more distributional information and, therefore, effectively mitigate the impact of future attacker's responses). 


In conclusion, we note that the constraints (\ref{cons: auxiliary probability}) and (\ref{cons: auxiliary expectation}) are necessary only if they refine the initial distributional constraints~in~$\mathcal{Q}$. A possible verification of this property is discussed in more detail within Section \ref{subsec: multi-stage MIP}. In the following example, we provide some intuition behind the auxiliary distributional constraints (\ref{eq: list of additional constraints}), as well as the first and the second research questions, \textbf{Q1} and \textbf{Q2}, from Section \ref{subsec: approach and contributions}. 
\begin{example} \label{ex: comparison} \upshape
	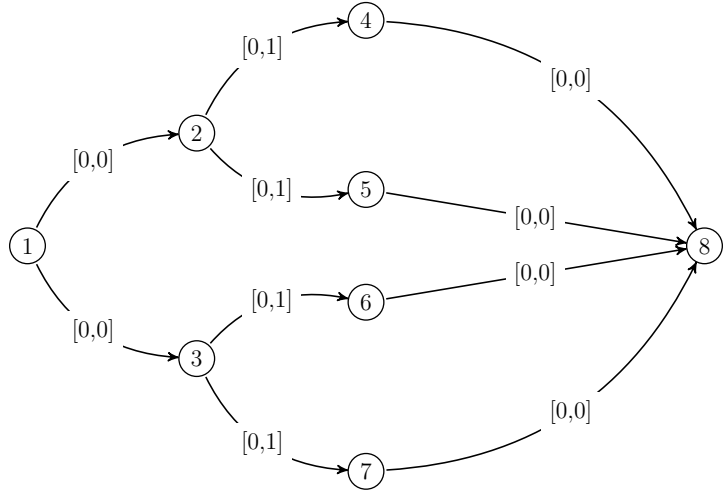
\begin{figure}
		\centering
		\begin{tikzpicture}[scale=0.75,transform shape]
		\Vertex[x=0,y=0]{1}
		\Vertex[x=3,y=2]{2}
		\Vertex[x=3,y=-2]{3}
		\Vertex[x=6,y=4]{4}
		\Vertex[x=6,y=1]{5}
		\Vertex[x=6,y=-1]{6}
		\Vertex[x=6,y=-4]{7}
		\Vertex[x=12,y=0]{8}
		\tikzstyle{EdgeStyle}=[post]
		\Edge[label=$[0\mbox{,}0]$](5)(8)
		\Edge[label=$[0\mbox{,}0]$](6)(8)
		\tikzstyle{EdgeStyle}=[post, bend left]
		\Edge[label=$[0\mbox{,}0]$](1)(2)
		\Edge[label=$[0\mbox{,}0]$](4)(8)
		\Edge[label=$[0\mbox{,}1]$](2)(4)
		\Edge[label=$[0\mbox{,}1]$](3)(6)
		\tikzstyle{EdgeStyle}=[post, bend right]
		\Edge[label=$[0\mbox{,}0]$](1)(3)
		\Edge[label=$[0\mbox{,}0]$](7)(8)
		\Edge[label=$[0\mbox{,}1]$](2)(5)
		\Edge[label=$[0\mbox{,}1]$](3)(7)
		\end{tikzpicture}
		\caption{\footnotesize The network used in Example \ref{ex: comparison}. The cost range is depicted inside each arc.}
		\label{fig: comparison}
	\end{figure} 
	We consider a network $G$ depicted in Figure \ref{fig: comparison} with $s = 1$ and $f = 8$. The arc costs are supposed to satisfy support constraints outlined inside each arc, i.e., $c_a \in [l_a, u_a]$ with probability~$1$. Furthermore, we introduce a unique linear expectation constraint given by:
	\begin{equation} \label{ex: eq: linear expectation constraint}
	\mathbb{E}_\mathbb{Q}\{\sum_{a \in A'} c_a\} \leq 1,
	\end{equation}
	where $A' := \{(2,4),(2,5),(3,6),(3,7)\}$. Next, we consider three different cases with respect to the additional distributional information \textit{observed by the user}: 
	\begin{itemize}
		\item \textbf{Case 1 (no additional information).} If the user picks a path \textit{before} the realization of $\mathbb{Q} \in \mathcal{Q}$, then the user's path contains exactly one arc from $A'$ and the attacker may set the expected of this arc equal to $1$.
		Then the worst-case expected loss incurred by the user is also equal to $1$.
		\item \textbf{Case 2 (full information).}
		 If the user picks a path \textit{after} the realization of $\mathbb{Q} \in \mathcal{Q}$, then the optimal attacker's decision is to set expected costs of $0.25$ for each arc in $A'$. Then the worst-case expected loss of the user equals $0.25$ (this value can be also viewed as some \textit{lower bound} on the user's expected~loss). 
		\item \textbf{Case 2 (successive revelation of information).} Assume that, if the user stays at node $2$, then it may compare the expected costs of $(2,4)$ and $(2,5)$. In other words, the attacker indicates whether the following linear expectation constraint holds or not: 
		\begin{equation} \nonumber
		\mathbb{E}_{\mathbb{Q}} \{c_{(2,4)} - c_{(2,5)}\} \leq 0.
		\end{equation}
		In this case it is optimal for the user to traverse through $(1,2)$ and then to pick an arc with the smallest expected cost. In fact, the expected cost of the latter arc cannot exceed $0.5$ due to the initial constraint (\ref{ex: eq: linear expectation constraint}), which results in the worst-case expected loss of $0.5$. 
	\end{itemize}	
	
	Comparing the first and the third cases we conclude that the user may benefit from using dynamic decisions. Furthermore, in view of Case 2, the user's expected loss in our setting can be bounded from below by an optimal objective function value of the \textit{max-min problem}, i.e., the one-stage problem (\ref{intro: single-stage DRO}) where the order of ``max'' and ``min'' operators is reversed. \vspace{-9mm}\flushright$\square$
\end{example}

\subsection{Multi-stage formulation} \label{subsec: multi-stage}

Following our discussion in Section \ref{subsec: related literature}, we encode all feasible attacker's responses to the constraints in $\mathcal{L}$ with binary vectors $\mathbf{r}_j \in \{0, 1\}^{|\mathcal{L}|}$ for $j \in \{1, \ldots, 2^{|\mathcal{L}|}\}$. In particular, for any fixed $j$ we observe that $r_{j m} = 1$, if the $m$-th constraint in $\mathcal{L}$ is satisfied and $r_{j m} = 0$, otherwise. For each vector of attacker's responses $\mathbf{r}_{j}$, $j \in \{1, \ldots, 2^{|\mathcal{L}|}\}$, we construct a new ambiguity set $\mathcal{Q}_{j}$ that combines the initial distributional constraints in $\mathcal{Q}$ with the auxiliary constraints induced by $\mathbf{r}_{j}$. For ease of exposition we refer to ``$v \geq w$'' as an opposite constraint to ``$v \leq w$''. 
Formally, for each $j \in \{1, \ldots, 2^{|\mathcal{L}|}\}$ we have:
\begin{equation} \nonumber
\mathcal{Q}_{j} := \mathcal{Q} \cap \Big\{\mathbb{Q} \in \mathcal{Q}_0(\mathbb{R}^{|A|}) \mbox{ s.t. } \begin{cases} \mbox{ the } m \mbox{-th constraint in } \mathcal{L}, \mbox{ if } r_{j m} = 1 \\ 
\mbox{ an opposite of the } m \mbox{-th constraint in } \mathcal{L}, \mbox{ if } r_{j m} = 0 \end{cases} \Big\}.
\end{equation}

Then, irrespective of the order of auxiliary distributional constraints in $\mathcal{L}$, the multi-stage DRSPP can be formulated as follows:
\begin{equation} \tag{\textbf{F}$_{ms}$} \label{multi-stage DRSPP}
\begin{gathered} 
	\max_{j \in \{1, \ldots, 2^{|\mathcal{L}|}\}} \min_{\mathbf{y}_j \in Y} \max_{\mathbb{Q}_j \in \mathcal{Q}_j} \mathbb{E}_{\mathbb{Q}_j } \{\mathbf{c}^\top \mathbf{y}_j\} \qquad \qquad \qquad \qquad \qquad \qquad \qquad \qquad \quad \\
	\mbox{s.t. } \mbox{non-anticipativity constraints with respect to } \mathbf{y}_j, j \in \{1, \ldots, 2^{|\mathcal{L}|}\},	
\end{gathered}
\end{equation}
\looseness-1 where $\mathbf{y}_j$, $j \in \{1, \ldots, 2^{|\mathcal{L}|}\}$, denotes a decision of the user under complete knowledge of the vector of attacker's responses $\mathbf{r}_j$  or, equivalently, the associated ambiguity set $\mathcal{Q}_j$. Formally, in (\ref{multi-stage DRSPP}) the user attempts to minimize its worst-case expected loss taking into account all possible realizations of attacker's responses and subject to some non-aniticipativity constraints. 
We need to enforce non-anticipativity as long as the user learns the attacker's responses associated with some node $i \in N \setminus \{f\}$ only when it reaches node $i$. In other words, for any fixed $j, \ell \in \{1, \ldots, 2^{|\mathcal{L}|}\}$ the user's paths $P_j$ and $P_\ell$ (induced by the incidence vectors $\mathbf{y}_j$ and $\mathbf{y}_\ell$, respectively) must coincide whenever the user is not able to distinguish between the ambiguity sets $\mathcal{Q}_j$ and $\mathcal{Q}_{\ell}$. 

 First, we construct non-anticipativity constraints for acyclic graphs. In this regard, for any fixed $j, \ell \in \{1, \ldots, 2^{|\mathcal{L}|}\}$, $j \neq \ell$, we denote by $N_{j, \ell} \subseteq N$ a set of nodes at which the user may learn the actual ambiguity set, either $\mathcal{Q}_j$ or $\mathcal{Q}_{\ell}$, that is enforced by the attacker. More specifically, for every node in $N_{j, \ell}$ there exists an associated distributional constraint in $\mathcal{L}$ such that the attacker's responses to this constraint in $\mathbf{r}_{j}$ and $\mathbf{r}_{\ell}$ are different. As a remark, in Example \ref{ex: comparison} we have $|\mathcal{L}| = 1$ and $N_{1,2} = \{2\}$. The following result holds.

\begin{proposition} \label{propostion: non-anticipativity acyclic} 
Assume that $j, \ell \in \{1, \ldots, 2^{|\mathcal{L}|}\}$, $j \neq \ell$, and let $G$ be a directed acyclic graph. Then for each node $i \in N \setminus N_{j, \ell}$ the following constraints ensure non-anticipativity:
\begin{subequations} \label{cons: non-anticipativity acyclic}
	\begin{eqnarray}
	y_{j,a} = y_{\ell,a} \quad \forall a \in FS_i,\mbox{ if all nodes in } N_{j, \ell} \mbox{ are reachable from } i, \label{cons: non-anticipativity 1 acyclic}\\
	\begin{gathered} \begin{rcases}
	|y_{j,a} - y_{\ell,a}| \leq \sum_{n \in \widetilde{N}_{j, \ell}}\sum_{a' \in FS_{n}} y_{j,a'} \quad \forall a \in FS_i,\\
	\; \; \qquad \mbox{if a subset of nodes } \widetilde{N}_{j, \ell} \subseteq N_{j, \ell} \mbox{ is not reachable from } i. \label{cons: non-anticipativity 2 acyclic} \end{rcases}
	\end{gathered}
	\end{eqnarray}
\end{subequations}
\begin{proof}
 We consider three particular situations:
\begin{enumerate}
	\item[(\textit{i})] Assume that for some node $i \in N \setminus N_{j, \ell}$ all nodes in $N_{j, \ell}$ are reachable from $i$. Then any node in $N_{j, \ell}$ can be visited by the user only after $i$ since the graph $G$ is acyclic. Hence, at node $i$ we have a lack of information to distinguish between the ambiguity sets $\mathcal{Q}_j$ and $\mathcal{Q}_\ell$, which results in the constraints (\ref{cons: non-anticipativity 1 acyclic}).
	\item[(\textit{ii})] Assume that for some node $i \in N \setminus N_{j, \ell}$ a subset of nodes $\widetilde{N}_{j, \ell} \subseteq N_{j, \ell}$ is not reachable from~$i$. In view of the first part of non-anticipativity constraints (\ref{cons: non-anticipativity 1 acyclic}) and the definition of $\widetilde{N}_{j, \ell}$, it is rather straightforward to check that either both paths induced by $\mathbf{y}_j$ and $\mathbf{y}_\ell$ contain \textit{some} node in $\widetilde{N}_{j,\ell}$ or both not. In fact, non-anticipativity at node $i$ is needed only if the user's path induced by $\mathbf{y}_j$, say $P_j$, does not contain any node in $\widetilde{N}_{j, \ell}$, i.e., $$\sum_{n \in \widetilde{N}_{j, \ell}}\sum_{a' \in FS_{n}} y_{j,a'} = 0,$$
	see constraints (\ref{cons: non-anticipativity 2 acyclic}). Otherwise, there exists some node $n \in P_j \cap \widetilde{N}_{j, \ell}$ and we have either $i \notin P_j$ or $i \in P_j$ and $i$ is visited after $n$; recall the definition of $\widetilde{N}_{j, \ell}$. In both cases the constraints (\ref{cons: non-anticipativity 1 acyclic}) and the flow conservation constraints (\ref{eq: flow conservation}) are sufficient to guarantee non-anticipativity and, thus, we make the constraints (\ref{cons: non-anticipativity 2 acyclic}) non-binding.
	\item[(\textit{iii})] If the user stays at some node $i \in N_{j, \ell}$, then it may distinguish between the ambiguity sets $\mathcal{Q}_j$ and $\mathcal{Q}_\ell$. Hence, we do not need to enforce non-anticipativity constraints.
	
\end{enumerate}
These observations conclude the proof. 
\end{proof} 
\end{proposition} 

We note that reachability in the formulation of Proposition \ref{propostion: non-anticipativity acyclic} can be checked efficiently for each pair of nodes, e.g., by Floyd-Warshall algorithm; see \cite{Cormen2009}. Furthermore, the constraints (\ref{cons: non-anticipativity 1 acyclic}) and (\ref{cons: non-anticipativity 2 acyclic}) are, in fact, linear constraints with respect to the decision variables.

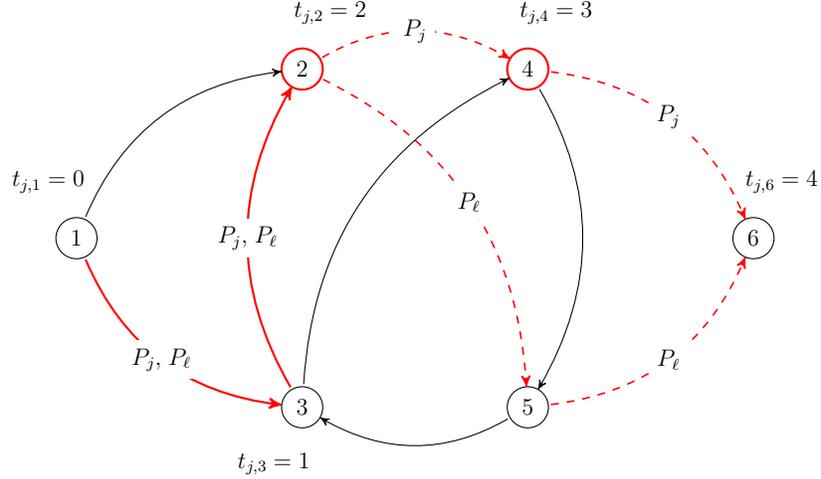
\begin{figure}
	
	\centering
	\begin{tikzpicture}[scale=0.75,transform shape]
	\tikzstyle{VertexStyle}=[fill=white,sloped]
	\Vertex[x=-0.5,y=1]{$t_{j,1} = 0$}
	\Vertex[x=3.5,y=-4]{$t_{j,3} = 1$}
	\Vertex[x=4.5,y=4]{$t_{j,2} = 2$}
	\Vertex[x=8.5,y=4]{$t_{j,4} = 3$}
	\Vertex[x=12.5,y=1]{$t_{j,6} = 4$}
	\tikzstyle{VertexStyle}=[draw = black, circle]
	\Vertex[x=0,y=0]{1}
	\Vertex[x=4,y=-3]{3}
	\Vertex[x=12,y=0]{6}
	\Vertex[x=8,y=-3]{5}
	\tikzstyle{VertexStyle}=[draw = red, circle, thick]
	\Vertex[x=4,y=3]{2}
	\Vertex[x=8,y=3]{4}
	\tikzstyle{EdgeStyle}=[post, bend left, thin]
	\Edge(1)(2)
	\Edge(5)(3)
	\Edge(3)(4)
	\Edge(4)(5)
	\tikzstyle{EdgeStyle}=[post, bend right, thick, red]
	\Edge[label=$P_j\mbox{,} \;P_{\ell}$](1)(3)
	\tikzstyle{EdgeStyle}=[post, bend left, thick, red]
	\Edge[label=$P_j\mbox{,} \;P_{\ell}$](3)(2)
	\tikzstyle{EdgeStyle}=[post, bend left, dashed, red]
	\Edge[label=$P_j$](2)(4)
	\Edge[label=$P_j$](4)(6)
	\tikzstyle{EdgeStyle}=[post, bend left, dashed, red]
	\Edge[label=$P_{\ell}$](2)(5)
	\tikzstyle{EdgeStyle}=[post, bend right, dashed, red]
	\Edge[label=$P_{\ell}$](5)(6)
	\end{tikzpicture}
	\caption{\footnotesize A pair of user's decisions satisfying the non-anticipativity requirement for some $j,\ell \in \{1, \ldots, 2^{|\mathcal{L}|}\}$. The set $N_{j, \ell}$ is assumed to contain two nodes, $2$ and $4$, depicted in red. Some feasible values of $\mathbf{t}_j$ are depicted outside the nodes of $P_j$.} 
	\label{fig: non-anticipativity}
\end{figure}

Unfortunately, the non-anticipativity constraints become substantially more complicated for general graphs containing directed cycles. For example, there may exist a cycle containing some node $i \in N \setminus N_{j,\ell}$ and $i' \in N_{j,\ell}$ for some fixed $j, \ell \in \{1, \ldots, 2^{|\mathcal{L}|}\}$. In this case we need to enforce non-anticipativity at node $i$ only if node $i'$ is visited by the user after $i$. For example, in Figure \ref{fig: non-anticipativity} we consider a general graph, where the non-anticipativity requirement (\ref{cons: non-anticipativity 1 acyclic}) for acyclic graphs is violated. Indeed, all nodes in $N_{j,\ell} = \{2, 4\}$ are reachable from node $5$, but we do not need to enforce non-anticipativity at this node, as node $2$ is visited before. 

In order to resolve the outlined issue, we introduce new variables $\mathbf{t}_j \in \mathbb{R}^{|N|}$,~$j~\in~\{1, \ldots, 2^{|\mathcal{L}|}\}$, related to a sequence at which nodes are visited by the user under scenario $j$; see constraints (\ref{cons: non-anticipativity time 1a})-(\ref{cons: non-anticipativity time 2}) used in the formulation of Proposition \ref{propostion: non-anticipativity general} below. For instance, in Figure \ref{fig: non-anticipativity} some feasible values of $\mathbf{t}_j$ associated with the user's path $P_j$ are depicted. As a remark, similar variables and constraints are used in the Miller–Tucker–Zemlin (MTZ) formulation of the travelling salesmen problem \cite{Miller1960}. We formulate the following result for general graphs. 

\begin{proposition} \label{propostion: non-anticipativity general} 
	Assume that $j,\ell \in \{1, \ldots, 2^{|\mathcal{L}|}\}$, $j \neq \ell$, and the graph $G$ is general. Then for each node $i \in N \setminus N_{j,\ell}$ the following constraints ensure non-anticipativity:
	\begin{subequations} \allowdisplaybreaks \label{cons: non-anticipativity general}
		\begin{align}
		& t_{j,s} = 0 \label{cons: non-anticipativity time 1a}\\
		& 0 \leq t_{j,i} \leq |N| - 1 \quad \forall i \in N \label{cons: non-anticipativity time 1b}\\
		& t_{j,a_1} - t_{j,a_2} \leq -1 + |N|(1 - y_{j,a}) \quad \forall a \in A \label{cons: non-anticipativity time 2} \\
		&\begin{rcases} |y_{j,a} - y_{\ell,a}| \leq \sum_{n \in N_{j,\ell}} \min\Big\{\max\{t_{j,i} - t_{j,n}; 0\} + 2 - \sum_{a' \in FS_{n}} y_{j,a'} - \qquad \\ 
		\qquad \qquad \qquad \qquad \qquad \qquad \qquad \qquad \qquad - \sum_{a' \in FS_{i}} y_{j,a'}; \sum_{a' \in FS_{n}} y_{j,a'}\Big\} \end{rcases}
		\quad \forall a \in FS_i. \label{cons: non-anticipativity 3}
		\end{align}
	\end{subequations}
	\begin{proof}
		 Henceforth, we fix $j,\ell \in \{1, \ldots, 2^{|\mathcal{L}|}\}$, $j \neq \ell$. First, we note that the decision vector $\mathbf{y}_j \in Y$ induces some simple $s-f$ path 
		$$P_j = \{n_1 \rightarrow n_2 \rightarrow \ldots \rightarrow n_{T+1}\},$$
	where $n_1 = s$, $n_{T+1} = f$ and any intermediate node is visited at most once. In fact, the constraints (\ref{cons: non-anticipativity time 1a})$-$(\ref{cons: non-anticipativity time 2}) assign labels $\mathbf{t}_j$ to each node $i \in N$ so that the label at $s$ is equal to $0$ and the labels at each subsequent node in $P_j$ are incremented at least by one compared to the predecessor node. Each label is also bounded from above by $|N| - 1$, which is the maximal possible length of a simple $s-f$ path in $G$. 

	Next, we note that the user's paths $P_j$ and $P_\ell$ (induced by the incidence vectors $\mathbf{y}_j$ and $\mathbf{y}_\ell$, respectively) must coincide until some node $n \in N_{j,\ell}$ is reached. This observation implies that non-anticipativity constraints are necessary only for the nodes of $P_j$ and $P_\ell$ that are visited before $n$. We consider two particular cases:
	\begin{enumerate}
		\item[(\textit{i})] Assume that the user stays at some node $i \in N \setminus N_{j,\ell}$. Then non-anticipativity at $i$ is needed only if \textit{all} nodes in $N_{j,\ell}$ are either not contained in $P_j$ or are visited after node $i$ assuming that $i \in P_j$. Hence, for each node $n \in N_{j,\ell}$ we verify whether this node is contained in $P_j$ and whether it is visited before or after node $i$. 
		
		Formally, we calculate the following indicators:
		$$w_{j,i,n} = \min\Big\{\max\{t_{j,i} - t_{j,n}; 0\} + 2 - \sum_{a' \in FS_{n}} y_{j,a'} - \sum_{a' \in FS_{i}} y_{j,a'}; \sum_{a' \in FS_{n}} y_{j,a'}\Big\}.$$
		We observe that, if $n$ is not contained in $P_j$, then $\sum_{a' \in FS_{n}} y_{j,a'} = 0$ and, hence, $w_{j,i,n} = 0$; recall that $\mathbf{y}_j$ must satisfy (\ref{eq: exlude cycles}). Otherwise, if $n$ is contained in $P_j$, then 
		$$w_{j,i,n} = \min\Big\{\max\{t_{j,i} - t_{j,n}; 0\} + 1 - \sum_{a' \in FS_{i}} y_{j,a'}; 1\Big\}.$$
		In this case, $w_{j,i,n} = 0$ only if $i \in P_j$, i.e., $\sum_{a' \in FS_{i}} y_{j,a'} = 1$, and $n$ is visited after $i$, i.e., $t_{j,i} \leq t_{j,n} - 1$. Otherwise, if $n$ is visited before $i$ or $i \notin P_j$, then $w_{j,i,n} = 1$. As a result, we conclude that non-anticipativity is guaranteed by the constraints (\ref{cons: non-anticipativity 3}). 
		\item[(\textit{ii})] If the user stays at some node $i \in N_{j,\ell}$, then it may distinguish between the ambiguity sets $\mathcal{Q}_j$ and $\mathcal{Q}_\ell$. Hence, we do not need to enforce non-anticipativity constraints.
	\end{enumerate}
	These observations conclude the proof. 
 \end{proof} 
\end{proposition} 

Similar to Proposition \ref{propostion: non-anticipativity acyclic} the non-anticipativity constraints (\ref{cons: non-anticipativity time 1a})$-$(\ref{cons: non-anticipativity 3}) can be reformulated as linear constraints but after applying some standard linearization techniques. We discuss this point as well as a linear MIP reformulation of (\ref{multi-stage DRSPP}) in Section \ref{sec: solution techniques}. In the remainder of this section we explore how the auxiliary distributional constraints given by (\ref{eq: list of additional constraints}) can be constructed from real data observations. 

\subsection{Construction of the auxiliary constraints from data} \label{subsec: constraints form data} 
The results of this section are motivated and similar to the related results of Section 2.2 in~\cite{Ketkov2021}, where the initial distributional constraints in $\mathcal{Q}$ are constructed from data. In the following, we assume that the attacker fixes some nominal distribution $\mathbb{Q}^* \in \mathcal{Q}$ and has access to a data set given by $\widehat{n} \in \mathbb{Z}_{>0}$ independent observations obtained from this distribution, i.e.,
\begin{equation} 
\widehat{\mathbf{C}} := \Big\{(\hat{c}_1^{(k)}, \ldots, \hat{c}_{|A|}^{(k)})^\top, \; k \in \{1, \ldots, \widehat{n}\} \Big\} \label{eq: data set}
\end{equation}
\looseness-1 We also note that an optimal solution of (\ref{multi-stage DRSPP}) provides a user's decision for \textit{every} possible sequence of attacker's responses. 

 A natural question arising is how the user can identify the actual sequence of attacker's responses. In this regard, \textit{after} solving the multi-stage formulation (\ref{multi-stage DRSPP}) we exploit the data set~(\ref{eq: data set}) to verify the auxiliary constraints~(\ref{eq: list of additional constraints}) dynamically at the respective nodes of the user's path. This procedure is referred to as a \textit{constraint verification procedure}.


\textbf{Verification of probability constraints.}
Suppose that the user stays at some node $i \in N \setminus \{f\}$, and needs to verify a probability constraint of the form 
\begin{equation} \label{eq: probability constraint to check} 
q^*_a := \mathbb{Q}^*_a\{c_a \in [\widetilde{l}_a, \widetilde{u}_a]\} \leq \widetilde{q}_a
\end{equation}
for some $a \in FS_i$. From Hoeffding inequality \cite{Hoeffding1994} we observe that for any $\varepsilon > 0$
\begin{equation} \label{eq: hoeffding 1}
\Pr\Big\{|q^*_a - \frac{1}{\widehat{n}}\sum_{k = 1}^{\widehat{n}} \chi_{a,k}| \geq \varepsilon\Big\} \leq 2 \exp(- 2\widehat{n} \varepsilon^2),
\end{equation}
 where $\frac{1}{\widehat{n}}\sum_{k = 1}^{\widehat{n}} \chi_{a,k}$ refers to an empirical probability that $c_a \in [\widetilde{l}_a, \widetilde{u}_a]$, i.e., 
\begin{equation} \nonumber
\chi_{a,k} = \begin{cases}
1, \mbox{ if } \hat{c}^{(k)}_a \in [\widetilde{l}_a, \widetilde{u}_a] \\
0, \mbox{ otherwise},
\end{cases}
\end{equation}
for each $k \in \{1, \ldots, \widehat{n}\}$.
Hence, with any prescribed confidence level $\gamma \in (0, 1)$ we may guarantee that 
$$q^*_a \in [\frac{1}{\widehat{n}}\sum_{k = 1}^{\widehat{n}} \chi_{a,k} - \varepsilon, \frac{1}{\widehat{n}}\sum_{k = 1}^{\widehat{n}} \chi_{a,k} + \varepsilon],$$
where the parameter $\varepsilon$ is defined by setting the right-hand side of (\ref{eq: hoeffding 1}) equal to $1 -\gamma$. As a remark, using the same arguments one may construct the initial probability constraints~(\ref{probability constraints}). 


Furthermore, if the number of samples, $\widehat{n}$, is sufficiently large (and, therefore, $\varepsilon$ is sufficiently small), then we may distinguish between the following two alternatives:
$$\frac{1}{\widehat{n}}\sum_{k = 1}^{\widehat{n}} \chi_{a,k} + \varepsilon \leq \widetilde{q}_a \; \mbox{ or } \; \frac{1}{\widehat{n}}\sum_{k = 1}^{\widehat{n}} \chi_{a,k} - \varepsilon > \widetilde{q}_a.$$
In the former case the constraint (\ref{eq: probability constraint to check}) holds with probability of at least $\gamma$; in the latter case it is violated with probability of at least $\gamma$. Importantly, in order to verify a constraint of the form (\ref{eq: probability constraint to check}), it suffices to know only the \textit{number} of random samples that belong to the interval~$[\widetilde{l}_a, \widetilde{u}_a]$. Also, if the constraint cannot be verified with $\widehat{n}$ samples, then one may generate the attacker's response at random or assume the worst-case scenario for the user. We exploit some variation of this idea within our numerical study in Section \ref{sec: comp study}. 


\textbf{Verification of linear expectation constraints.}
In a similar way, one may verify a linear expectation constraint of the form:
\begin{equation} \label{eq: expectation constraint to check}
\mathbb{E}_{{\mathbb{Q}}^*}\{\sum_{a \in FS_i} p_a c_a\} \leq p_0
\end{equation}
for some $p_a \in \mathbb{R}$, $a \in FS_i$, and $p_0 \in \mathbb{R}$. 
That is, by Hoeffding inequality we observe that:
\begin{equation} \label{eq: hoeffding 2}
\Pr\Big\{|\mathbb{E}_{{\mathbb{Q}}^*}\{\sum_{a \in FS_i} p_a c_a\} - \frac{1}{\widehat{n}}\sum_{k = 1}^{\widehat{n}} \sum_{a \in FS_i} p_a \hat{c}^{(k)}_a | \geq \varepsilon\Big\} \leq 2 \exp\Big(-2\widehat{n} \Big(\frac{\varepsilon}{\overline{S}_i - \underline{S}_i}\Big)^2\Big),
\end{equation}
where $\overline{S}_i = \sum_{a \in FS_i} \max_{c_a \in [l_a, u_a]} p_a c_a$, $\underline{S}_i = \sum_{a \in FS_i} \min_{c_a \in [l_a, u_a]} p_a c_a$ and $i \in N \setminus \{f\}$. Then for any confidence level $\gamma \in (0, 1)$ we may guarantee that:
$$\mathbb{E}_{{\mathbb{Q}}^*}\{\sum_{a \in FS_i} p_a c_a\} \in [\frac{1}{\widehat{n}}\sum_{k = 1}^{\widehat{n}} \sum_{a \in FS_i} p_a \hat{c}^{(k)}_a - \varepsilon, \frac{1}{\widehat{n}}\sum_{k = 1}^{\widehat{n}} \sum_{a \in FS_i} p_a \hat{c}^{(k)}_a + \varepsilon],$$
where $\varepsilon$ is defined by setting the right-hand side of (\ref{eq: hoeffding 2}) equal to $1 - \gamma$. As in the previous case, the above procedure can be slightly modified to construct the initial expectation constraints~(\ref{linear expectation constraints}).

We note that the constraint (\ref{eq: expectation constraint to check}) can be verified following the same arguments as for the probability constraint (\ref{eq: probability constraint to check}). However, in contrast to (\ref{eq: hoeffding 1}), the inequality in (\ref{eq: hoeffding 2}) exploits \textit{linear combinations} of arc costs with respect to the arcs $a \in FS_i$. In the following example we demonstrate how the auxiliary constraints (\ref{cons: auxiliary probability}) and (\ref{cons: auxiliary expectation}) can be constructed and verified by using information from Bluetooth sensors. 

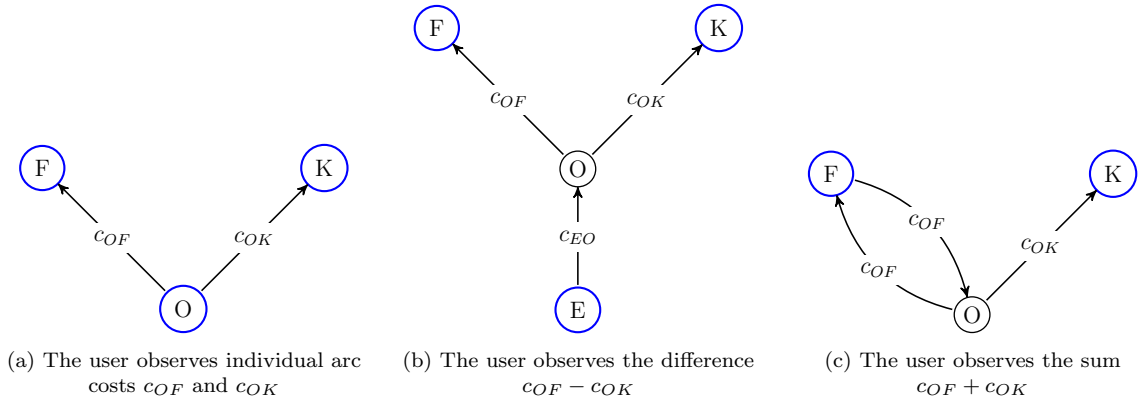
\begin{figure}
	\begin{subfigure}[b]{0.3\textwidth}
		\centering
		\captionsetup{justification=centering}
		\begin{tikzpicture}[scale=0.75,transform shape]
		\tikzstyle{VertexStyle}=[draw = blue, circle, thick]
		\Vertex[x=0,y=0]{O}
		\Vertex[x=-2.5,y=2.5]{F}
		\Vertex[x=2.5,y=2.5]{K}
		\tikzstyle{EdgeStyle}=[post]
		\Edge[label=$c_{OF}$](O)(F)
		\Edge[label=$c_{OK}$](O)(K)
		\end{tikzpicture}
		\caption{The user observes individual arc costs $c_{OF}$ and $c_{OK}$}
		\label{fig: ex2a}
	\end{subfigure}
	\begin{subfigure}[b]{0.3\textwidth}
		\centering
		\captionsetup{justification=centering}
		\begin{tikzpicture}[scale=0.75,transform shape]
		\Vertex[x=0,y=0]{O}
		\tikzstyle{VertexStyle}=[draw = blue, circle, thick]
		\Vertex[x=0,y=-2.5]{E}
		\Vertex[x=-2.5,y=2.5]{F}
		\Vertex[x=2.5,y=2.5]{K}
		\tikzstyle{EdgeStyle}=[post]
		\Edge[label=$c_{OF}$](O)(F)
		\Edge[label=$c_{OK}$](O)(K)
		\Edge[label=$c_{EO}$](E)(O)
		\end{tikzpicture} 		
		\caption{The user observes the difference $c_{OF} - c_{OK}$}
		\label{fig: ex2b}
	\end{subfigure}
	\begin{subfigure}[b]{0.3\textwidth}
		\centering
		\captionsetup{justification=centering}
		\begin{tikzpicture}[scale=0.75,transform shape]
		\Vertex[x=0,y=0]{O}
		\tikzstyle{VertexStyle}=[draw = blue, circle, thick]
		\Vertex[x=-2.5,y=2.5]{F}
		\Vertex[x=2.5,y=2.5]{K}
		\tikzstyle{EdgeStyle}=[post]
		\Edge[label=$c_{OK}$](O)(K)
		\tikzstyle{EdgeStyle}=[post, bend left]
		\Edge[label=$c_{OF}$](F)(O)
		\tikzstyle{EdgeStyle}=[post, bend left]
		\Edge[label=$c_{OF}$](O)(F)
		\end{tikzpicture}
		\caption{The user observes the sum $c_{OF} + c_{OK}$}
		\label{fig: ex2c}
	\end{subfigure}
	\caption{In each figure we consider different allocations of Bluetooth sensors, i.e., the sensors are placed at nodes highlighted in blue. The arc travel times $c_{OF}$ and $c_{FO}$ are assumed to be the same.}
	\label{fig: ex2}
\end{figure}

\begin{example} \label{example 2}
	\upshape	
	We envision a network, where the arc costs/travel times are measured using Bluetooth sensors; see, e.g., \cite{Asudegi2013, Haghani2010}. The sensors are placed at particular nodes of the network and one may quantify travel times between two successive sensors; see Figure \ref{fig: ex2}.  
	For example, in Figure \ref{fig: ex2a} the sensors are placed at nodes $O$, $F$ and $K$, which yields that the user may receive random observations of individual arc travel times $c_{OF}$ and~$c_{OK}$. On the other hand, in Figures \ref{fig: ex2b} and \ref{fig: ex2c}
	the user may only observe some linear combinations of individual arc travel times (we additionally assume that $c_{OF} = c_{FO}$ in Figure \ref{fig: ex2c}). 
	
	 Then, by leveraging the observed linear combinations of arc travel times and the constraint verification procedure, the user may verify linear expectation constraints of the form (\ref{cons: auxiliary expectation}). Analogously, if random observations of individual arc travel times are subject to interval uncertainty, then one may verify probability constraints of the form (\ref{cons: auxiliary probability}). \vspace{-9.5mm}\flushright$\square$ 
\end{example}

We conclude that the auxiliary distributional constraints (\ref{cons: auxiliary probability}) and (\ref{cons: auxiliary expectation}) can be motivated in a rather straightforward way by the outlined sensor-related application context. 
Some further intuition behind the choice of the auxiliary distributional constraints (\ref{cons: auxiliary probability}) and (\ref{cons: auxiliary expectation}) is implied by our solution procedure for the multi-stage problem (\ref{multi-stage DRSPP}) and is discussed within Section \ref{sec: solution techniques}.

\section{Solution approach} \label{sec: solution techniques}
This section is organized as follows. First, in Section \ref{subsec: one-stage MIP} we briefly describe the key theoretical results for the one-stage model in \cite{Ketkov2021}. Then in Section \ref{subsec: multi-stage MIP} we use the aforementioned results as well as Propositions \ref{propostion: non-anticipativity acyclic} and \ref{propostion: non-anticipativity general} to reformulate the multi-stage DRSPP (\ref{multi-stage DRSPP}) as a linear MIP problem. 

\subsection{MIP reformulation of the one-stage problem} \label{subsec: one-stage MIP}
 We recall that the one-stage DRSPP in \cite{Ketkov2021} can be expressed as: 
\begin{equation} \label{one-stage DRSPP} \tag{\textbf{F}$_{os}$}
z^*_{static} := \min_{\mathbf{y} \in Y} \max_{\mathbb{Q} \in \mathcal{Q}} \mathbb{E}_{\mathbb{Q}} \{\mathbf{c}^\top \mathbf{y}\},
\end{equation}
where the ambiguity set $\mathcal{Q}$ and the set of feasible decisions $Y$ are given by (\ref{ambiguity set}) and~(\ref{path flow constraints}), respectively. The optimization problem (\ref{one-stage DRSPP}) is also referred to as a \textit{static problem} and its optimal objective function value is denoted as $z^*_{static}$. 

It turns out that (\ref{one-stage DRSPP}) can be recast as a robust shortest path problem with some polyhedral uncertainty set. More specifically, the following result holds. 

\begin{theorem} \label{theorem 1} (Theorem 2 in \cite{Ketkov2021}) Assume that
	\begin{equation} \nonumber
	\mathcal{S}_0: = \{\overline{\mathbf{c}} \in \mathbb{R}^{|A|}:\mbox{ } \mathbf{L} \leq \overline{\mathbf{c}} \leq \mathbf{U}; \quad \mathbf{B}\overline{\mathbf{c}} \leq \mathbf{b}\},
	\end{equation}
where for each $a \in A$	
\begin{subequations}
	\begin{align}
	& L_a := \min_{\mathbb{Q}_a \in \widetilde{\mathcal{Q}}_a} \mathbb{E}_{\mathbb{Q}_a}\{c_a\}, \label{best case expected cost}\\
	& U_a := \max_{\mathbb{Q}_a \in \widetilde{\mathcal{Q}}_a} \mathbb{E}_{\mathbb{Q}_a}\{c_a\}, \label{worst case expected cost}
	\end{align}
\end{subequations}
and 
\begin{equation} \nonumber
\widetilde{\mathcal{Q}}_a := \Big\{\mathbb{Q}_a \in \mathcal{Q}_0(\mathbb{R}) \mbox{: }
\mathbb{Q}_a\{c_a \in [l^{(i)}_a, u^{(i)}_a]\} \in [\underline{q}^{(i)}_a, \overline{q}^{(i)}_a] \quad \forall i \in \mathcal{D}_a
\Big\}.
\end{equation}
Then the DRSPP of the form (\ref{intro: single-stage DRO}) is equivalent to the following robust shortest path problem with polyhedral uncertainty:
	\begin{equation} \label{robust formulation} 
	\min_{\mathbf{y} \in Y} \max_{\overline{\mathbf{c}} \in \mathcal{S}_0} \; \overline{\mathbf{c}}^\top \mathbf{y}.
	\end{equation}
\end{theorem}

The intuition behind Theorem \ref{theorem 1} can be explained as follows. Since the objective function in (\ref{one-stage DRSPP}) is linear, it can be seen as a function of expected costs, $\overline{\mathbf{c}} := \mathbb{E}_{\mathbb{Q}}\{\mathbf{c}\}$, i.e.,
$$\mathbb{E}_{\mathbb{Q}}\{\mathbf{c}^\top \mathbf{y}\} = (\mathbb{E}_{\mathbb{Q}}\{\mathbf{c}\})^\top \mathbf{y} = \overline{\mathbf{c}}^\top \mathbf{y}.$$
One may also show that the individual probability constraints (\ref{probability constraints}) for each $a \in A$ can be expressed in terms of box constraints 
$$L_a \leq \overline{c}_a \leq U_a$$
with respect to the expected costs $\overline{c}_a$, which yields the desired result. 

In order to resolve the moment problems (\ref{best case expected cost}) and (\ref{worst case expected cost}), we make the following additional assumptions:
\begin{itemize}
	\item[\textbf{A4.}] For each $a \in A$ there exists a marginal distribution $\mathbb{Q}_a \in \mathcal{Q}_0(\mathbb{R})$ such that
	$$\mathbb{Q}_a\{l^{(i)}_a \leq c_a \leq u^{(i)}_a\} \in (\underline{q}^{(i)}_a, \overline{q}^{(i)}_a),$$
	whenever $\underline{q}^{(i)}_a < \overline{q}^{(i)}_a$, $i \in \mathcal{D}_a$.
	\item[\textbf{A5.}] For each $a \in A$ and any pair of subintervals in (\ref{probability constraints}), namely, $[l^{(i_1)}_a, u^{(i_1)}_a]$ and $[l^{(i_2)}_a, u^{(i_2)}_a]$, $i_1, i_2 \in \mathcal{D}_a$, we have $l^{(i_1)}_a \neq u^{(i_2)}_a$ and $l^{(i_2)}_a \neq u^{(i_1)}_a$.
\end{itemize}
 Simply speaking, 
Assumptions \textbf{A4} and \textbf{A5} guarantee that strong duality for the moment problems (\ref{best case expected cost}) and (\ref{worst case expected cost}) holds. This prerequisite allows us to obtain finite linear programming reformulations of (\ref{best case expected cost}) and (\ref{worst case expected cost}); we refer to Lemma 2 in \cite{Ketkov2021} for further details. 

Summarizing the discussion above the one-stage formulation in \cite{Ketkov2021} can be tackled by solving $2|A|$ linear programming problems and a single robust shortest path problem (\ref{robust formulation}). Importantly, the latter problem can be recast as a linear MIP problem by dualizing the second-level linear programming problem. Formally, the following result holds.

\begin{proposition}[\textit{Theorem 3 in \cite{Ketkov2021}}] \label{proposition static MIP}
Let
\begin{equation} \label{polyhedral uncertainty set}
\mathcal{S}_0: = \Big\{\overline{\mathbf{c}} \in \mathbb{R}^{|A|}:\; \mathbf{L} \leq \overline{\mathbf{c}} \leq \mathbf{U}; \quad \mathbf{B}\overline{\mathbf{c}} \leq \mathbf{b}\Big\} = \Big\{\overline{\mathbf{c}} \in \mathbb{R}^{|A|}:\; \mathbf{B}_0\overline{\mathbf{c}} \leq \mathbf{b}_0\Big\}.
\end{equation}
Then the one-stage DRSPP (\ref{one-stage DRSPP}) admits the following mixed-integer programming reformulation:
\begin{equation} \label{static mip formulation} 
z^*_{static} = \min_{\mathbf{y}, \boldsymbol{\lambda}} \Big\{ \mathbf{b}_0^\top \boldsymbol{\lambda}: \; \boldsymbol{\lambda} \geq 0, \; -\mathbf{y} + \mathbf{B}_0^\top \boldsymbol{\lambda} = 0, \; \mathbf{y} \in Y \Big\}.
\end{equation}
\end{proposition}

\subsection{MIP reformulation of the multi-stage problem} \label{subsec: multi-stage MIP}
The key observation that we use next is that each family of distributions $\mathcal{Q}_j$, $j \in \{1, \ldots, 2^{|\mathcal{L}|}\}$ (corresponding to a vector of attacker's responses $\mathbf{r}_j \in \{0, 1\}^{|L|}$), contains the same types of distributional constraints as those in $\mathcal{Q}$; recall Assumption \textbf{A3}. Therefore, in view of Theorem~\ref{theorem 1}, each ambiguity set $\mathcal{Q}_j$, $j \in \{1, \ldots, 2^{\mathcal{L}}\}$, can be seen as some polyhedral uncertainty set 
\begin{equation} \label{partitioned polyhedrons}
\mathcal{S}_j := \Big\{\overline{\mathbf{c}} \in \mathbb{R}^{|A|}:\; \mathbf{B}_j \overline{\mathbf{c}} \leq \mathbf{b}_j \Big\} \subseteq \mathcal{S}_0
\end{equation}
in terms of expected costs. 
Next, we provide a linear MIP reformulation of the multi-stage problem~(\ref{multi-stage DRSPP}).

\begin{theorem} \label{theorem2}
Let $G$ be a general graph and assume that each ambiguity set $\mathcal{Q}_j$, $j \in \{1, \ldots, 2^{|\mathcal{L}|}\}$, is described by the polyhedral uncertainty set $\mathcal{S}_j$ given by (\ref{partitioned polyhedrons}).
Also, set~$M_1~=~|N|~-~1$. 
Then the multi-stage DRSPP (\ref{multi-stage DRSPP}) can be reformulated as the following linear MIP problem:
\begin{subequations} \allowdisplaybreaks \label{MIP reformulation} 
	\begin{align} 
	& z^*_{dynamic} = \min_{\mathbf{y}, \mathbf{t}, \mathbf{v}, \mathbf{w}, z} z \\ & \mbox{\upshape{s.t.} }
	\begin{rcases} t_{j,s} = 0 \\
	0 \leq t_{j,i} \leq |N| - 1 \quad \forall i \in N \\
	t_{j,a_1} - t_{j,a_2} \leq -1 + |N|(1 - y_{j,a}) \quad \forall a \in A \end{rcases} \quad \forall j \in \{1, \ldots, 2^{|\mathcal{L}|}\} \label{cons: times}\\ 
	& \begin{rcases}
	v_{j,i,n} \geq 0 \\
	v_{j,i,n} \geq t_{j,i} - t_{j,n} \\
	v_{j,i,n} \leq M_1 \tilde{v}_{j,i,n} \\
	v_{j,i,n} \leq t_{j,i} - t_{j,n} + M_1(1 - \tilde{v}_{j,i,n}) \\
	w_{j,i,n} \leq v_{j,i,n} + 2 - \sum_{a' \in FS_{n}} y_{j,a'} - \sum_{a' \in FS_{i}} y_{j,a'} \\
	w_{j,i,n} \leq \sum_{a' \in FS_{n}} y_{j,a'} \\
	-\sum_{\tilde{n} \in N_{j,\ell}} w_{j,i,\tilde{n}} \leq y_{j,a} - y_{\ell,a} \leq \sum_{\tilde{n} \in N_{j,\ell}} w_{j,i,\tilde{n}} \quad \forall a \in FS_i \\
	\tilde{v}_{j,i,n} \in \{0, 1\} \\
	\end{rcases} \quad \begin{array}{l} \forall j,\ell \in \{1, \ldots, 2^{|\mathcal{L}|}\},\\
	\forall n \in N_{j,\ell}, \; \\ \forall i \in N \setminus N_{j,\ell}, \\
	j \neq \ell \end{array} \label{cons: non-anticipativity general linearized} \\
	& \begin{rcases} z \geq \mathbf{b}_j^\top \boldsymbol{\lambda}_j \\
	-\mathbf{y}_j + \mathbf{B}_j^\top \boldsymbol{\lambda}_j = 0 \\
	\boldsymbol{\lambda}_j \geq 0 \\
	\mathbf{y}_j \in Y \end{rcases} \quad \forall j \in \{1, \ldots, 2^{|\mathcal{L}|}\}. \label{cons: worst-case expectation dualized}
	\end{align}
\end{subequations}
\begin{proof}
First, we provide an epigraph reformulation of (\ref{multi-stage DRSPP}) by using an auxiliary variable $z \in \mathbb{R}$, that is, 
\begin{subequations} \allowdisplaybreaks \label{multi-stage DRSPP epigraph}
	\begin{align}
	& \min z \\
	\mbox{s.t. } & \mbox{non-anticipativity constraints (\ref{cons: non-anticipativity general}),} \\
	& z \geq \min_{\mathbf{y}_j \in Y} \max_{\mathbb{Q}_j \in \mathcal{Q}_j} \mathbb{E}_{\mathbb{Q}_j} \{\mathbf{c}^\top \mathbf{y}_j\} \quad \forall j \in \{1, \ldots, 2^{|\mathcal{L}|}\}. \label{cons: epigraph 1}
	\end{align}
\end{subequations}
Next, the minimum in the right-hand side of (\ref{cons: epigraph 1}) indicates that for each $j \in \{1, \ldots, 2^{|\mathcal{L}|}\}$ there exists some $\widetilde{\mathbf{y}}_j \in Y$ such that 
$$z \geq \max_{\mathbb{Q}_j \in \mathcal{Q}_j} \mathbb{E}_{\mathbb{Q}_j} \{\mathbf{c}^\top \widetilde{\mathbf{y}}_j\}.$$
Hence, the path flow constraints $\mathbf{y}_j \in Y$ can be shifted to the first-level problem constraints by omitting the minimum in the right-hand side of (\ref{cons: epigraph 1}). In other words, the optimization problem (\ref{multi-stage DRSPP epigraph}) admits the following equivalent reformulation:
\begin{subequations} \label{multi-stage DRSPP epigraph 2} \allowdisplaybreaks
	\begin{align}
	& \min_{z, \mathbf{y}} z \\
	\mbox{s.t. } & \mbox{non-anticipativity constraints (\ref{cons: non-anticipativity general}), } \\
	& z \geq \max_{\mathbb{Q}_j \in \mathcal{Q}_j}\mathbb{E}_{\mathbb{Q}_j} \{\mathbf{c}^\top \mathbf{y}_j\} \quad \forall j \in \{1, \ldots, 2^{|\mathcal{L}|}\} \label{cons: epigraph 2, 1} \\
	& \mathbf{y}_j \in Y \quad \forall j \in \{1, \ldots, 2^{|\mathcal{L}|}\}.
	\end{align}
\end{subequations}

Then, in view of Theorem \ref{theorem 1}, the constraints (\ref{cons: epigraph 2, 1}) for each $j \in \{1, \ldots, 2^{|\mathcal{L}|}\}$ can be reformulated as polyhedral constraints with respect to expected costs $\overline{\mathbf{c}}_j = \mathbb{E}_{\mathbb{Q}_j} \{\mathbf{c}\}$, i.e.,
\begin{equation} \label{eq: robust reformulation for each family of distributions}
z \geq \max_{\overline{\mathbf{c}}_j \in \mathcal{S}_j} \; \overline{\mathbf{c}}_j^\top \mathbf{y}_j.
\end{equation}
Note that the maximization problem in the right-hand side of (\ref{eq: robust reformulation for each family of distributions}) is a linear program. 
Therefore, by strong duality we have:
\begin{equation} \nonumber
\max_{\overline{\mathbf{c}}_j \in \mathcal{S}_j} \; \overline{\mathbf{c}}_j^\top \mathbf{y}_j = \min_{\mathbf{y}_j, \boldsymbol{\lambda}_j} \Big\{ \mathbf{b}_j^\top \boldsymbol{\lambda}_j: \; \boldsymbol{\lambda}_j \geq 0, \; -\mathbf{y}_j + \mathbf{B}_j^\top \boldsymbol{\lambda}_j = 0 \Big\}.
\end{equation}
This observation and eliminating the minimum in the dual problem yields the following reformulation of (\ref{multi-stage DRSPP}):
\begin{subequations} \allowdisplaybreaks \label{MIP reformulation 0} 
\begin{align}
& \min_{\mathbf{y}, z} z\\
\mbox{\upshape{s.t.} } & \mbox{non-anticipativity constraints (\ref{cons: non-anticipativity general}), } \\
& \begin{rcases} z \geq \mathbf{b}_j^\top \boldsymbol{\lambda}_j \\
-\mathbf{y}_j + \mathbf{B}_j^\top \boldsymbol{\lambda}_j = 0 \\
\boldsymbol{\lambda}_j \geq 0 \\
\mathbf{y}_j \in Y \end{rcases} \quad \forall j \in \{1, \ldots, 2^{|\mathcal{L}|}\}. 
\end{align}
\end{subequations}

Finally, we linearize the non-anticipativity constraints (\ref{cons: non-anticipativity general}) defined in Proposition \ref{propostion: non-anticipativity general}. In this regard, for any fixed $j,\ell \in \{1, \ldots, 2^{|\mathcal{L}|}\}$, $n \in N_{j,\ell}$ and $i \in N \setminus N_{j,\ell}$ we introduce new variables $v_{j,i,n}$ and $w_{j,i,n}$ such that:
\begin{subequations}
 \begin{align}
 & v_{j,i,n} = \max\{0; t_{j,i} - t_{j,n}\} \label{eq: linearize max} \\ 
 & w_{j,i,n} = \min\{v_{j,i,n} + 2 - \sum_{a' \in FS_{n}} y_{j,a'} - \sum_{a' \in FS_{i}} y_{j,a'}; \sum_{a' \in FS_{n}} y_{j,a'} \}. \label{eq: linearize min} 
 \end{align}
\end{subequations}
The maximum in (\ref{eq: linearize max}) can be linearized by introducing new binary variables $\tilde{v}_{i,j,n} \in \{0, 1\}$ and the following linear constraints \cite{Nemhauser1988}:
\begin{align}
& v_{j,i,n} \geq 0 \nonumber \\
& v_{j,i,n} \geq t_{j,i} - t_{j,n} \nonumber \\
& v_{j,i,n} \leq M_1 \tilde{v}_{i,j,n} \nonumber \\
& v_{j,i,n} \leq t_{j,i} - t_{j,n} + M_1(1 - \tilde{v}_{j,i,n}). \nonumber
\end{align}
Specifically, as long as $|t_{j,i} - t_{j,n}| \leq |N| - 1$, we set $M_1 = |N| - 1$. 

Next, the constraint (\ref{cons: non-anticipativity 3}) implies that the sum of $w_{j,i,n}$ over $n \in N_{j, \ell}$ is bounded from below. Hence, in order to linearize the minimum in (\ref{eq: linearize min}) it is sufficient to enforce that
\begin{align}
& w_{j,i,n} \leq v_{j,i,n} + 2 - \sum_{a' \in FS_{n}} y_{j,a'} - \sum_{a' \in FS_{i}} y_{j,a'} \nonumber \\
& w_{j,i,n} \leq \sum_{a' \in FS_{n}} y_{j,a'}. \nonumber
\end{align}
As a result, the multi-stage DRSPP (\ref{multi-stage DRSPP}) can be expressed as a linear MIP problem (\ref{MIP reformulation}).
\end{proof}
\end{theorem}

Theorem \ref{theorem2} addresses the research question \textbf{Q3} formulated in Section \ref{subsec: approach and contributions}. Thus, the multi-stage DRSPP (\ref{multi-stage DRSPP}) can be solved at hand by using off-the-shelf MIP solvers even for graphs that contain cycles. In particular, a linear MIP reformulation for acyclic graphs can be derived from (\ref{MIP reformulation}) by using the non-anticipativity constraints (\ref{cons: non-anticipativity acyclic}) instead of the constraints (\ref{cons: times}) and (\ref{cons: non-anticipativity general linearized}). We point out that, if the number of auxiliary constraints, $|\mathcal{L}|$, is fixed, then the number of variables and constraints in (\ref{MIP reformulation}) is polynomial in the size of the network.

In addition, we observe that the auxiliary constraints (\ref{eq: list of additional constraints}) refine the initial distributional constraints in $\mathcal{Q}$ only if the resulting polyhedral uncertainty sets (\ref{partitioned polyhedrons}) are~non-empty. This property can be readily checked by solving the related linear feasibility problems for each polyhedron. Meanwhile, Theorem \ref{theorem2} remains valid even if some polyhedral uncertainty sets are empty. More precisely, if a set $\mathcal{S}_j$ is empty for some $j \in \{1, \ldots, 2^{|\mathcal{L}|}\}$, then $$\min_{\mathbf{y}_j, \boldsymbol{\lambda}_j} \Big\{ \mathbf{b}_j^\top \boldsymbol{\lambda}_j: \; \boldsymbol{\lambda}_j \geq 0, \; -\mathbf{y}_j + \mathbf{B}_j^\top \boldsymbol{\lambda}_j = 0 \Big\} = - \infty$$
by strong duality and, therefore, the associated constraints (\ref{cons: worst-case expectation dualized}) are non-binding. 

Finally, Theorem~\ref{theorem2} is also applicable to some other user's objective criteria, for which the constraints (\ref{cons: epigraph 2, 1}) admit a linear programming dual reformulation. In view of our results for the one-stage model in \cite{Ketkov2021}, the user may potentially optimize the conditional value at risk \cite{Rockafellar2000} or some other optimized certainty equivalent risk measures \cite{Hanasusanto2016, Wiesemann2014}. We leave these extensions as a possible direction of future research.

\section{Computational experiments} \label{sec: comp study}
In this section we attempt to answer the research questions \textbf{Q1} and \textbf{Q2}, i.e., we explore whether it is favorable for the user to employ \textit{dynamic} decisions and how much profit the user can gain by leveraging such decisions in place of \textit{static} decisions provided by the one-stage formulation (\ref{one-stage DRSPP}). In this regard, we analyze the role of the auxiliary distributional constraints (\ref{eq: list of additional constraints}) in relation to both the quality of dynamic decisions and tractability of the MIP reformulation~(\ref{MIP reformulation}).

The remainder of this section is organized as follows. In Section \ref{subsec: perfomance metrics} we introduce some performance metrics to assess the quality of dynamic decisions. Section \ref{subsec: test instances} describes our construction of test instances including the structure of the network and distributional constraints. Finally, in Section \ref{subsec: result for different types of networks} we provide our numerical results and the related discussion. 
All experiments are performed on a PC with \textit{CPU i5-7200U} and \textit{RAM 8 GB}. MIP problems are solved in Java with \textit{CPLEX 20.1}. 

\subsection{Performance metrics} \label{subsec: perfomance metrics}
 In view of our discussion in Section \ref{subsec: constraints form data}, we propose the following two-step procedure for validation of our approach. In the first step, we construct the initial ambiguity set $\mathcal{Q}$ and solve the MIP reformulation (\ref{MIP reformulation}) to obtain an optimal user's decision for every possible vector of attacker's responses. In the second step (which is referred to as a \textit{constraint verification~procedure}), we attempt to identify an actual sequence of attacker's responses by verifying the auxiliary distributional constraints~(\ref{eq: list of additional constraints}) dynamically at the respective nodes of the user's path.

 The quality of dynamic decisions in the first step is estimated using a relative gap of the form:
\begin{equation} \label{eq: relative gap 1} \tag{\textbf{G}$_1$}
\rho_1 := 100 \times \frac{z^*_{static} - z^*_{dynamic}}{z^*_{static} - z_{lower}}, 
\end{equation}
where $z_{lower}$ provides some lower bound on the optimal objective function value of the multi-stage problem (\ref{multi-stage DRSPP}). Following Example \ref{ex: comparison}, we select $z_{lower}$ as an optimal objective function value of related the max-min problem, i.e.,
\begin{equation} \label{max-min problem} 
z^*_{dynamic} \geq z_{lower} := \max_{\mathbb{Q} \in \mathcal{Q}} \min_{\mathbf{y} \in Y} \mathbb{E}_{\mathbb{Q}}\{\mathbf{c}^\top \mathbf{y}\}.
\end{equation}
Using the proof of Theorem \ref{theorem 1} it is quite easy to show that the max-min problem in (\ref{max-min problem}) admits a linear programming reformulation (the details are omitted for~brevity). 

The relative gap $\rho_1 \in [0, 100]$ quantifies the user's profit (in percentages) obtained from using dynamic decisions, in relation to the static problem formulation (\ref{one-stage DRSPP}). For instance, in Example \ref{ex: comparison} we have $z^*_{static} = 1$ and $z_{lower} =~ 0.25$. At the same time, it is rather straightforward to verify that $z^*_{dynamic} \geq 0.5$ irrespective of the auxiliary distributional constraints observed by the user. Hence, in this example $\rho_1 \leq 50\%$ and the value of 100$\%$ is not always achievable. 

On other hand, the quality of the constraint verification procedure can be defined as:
\begin{equation} \tag{\textbf{G}$_2$}
\rho_2 := 100 \times \frac{z^*_{dynamic} - \widetilde{z}_{dynamic}}{z^*_{static} - z_{lower}}, \label{eq: relative gap 2}
\end{equation}
where $\widetilde{z}_{dynamic}$ is our estimate of the worst-case expected loss \textit{after} the constraint verification procedure. More precisely, $\widetilde{z}_{dynamic}$ can be computed by solving the following linear programming problem:
\begin{equation} \nonumber
\widetilde{z}_{dynamic} := \max_{\mathbb{Q} \in \widetilde{\mathcal{Q}}} \mathbb{E}_{\mathbb{Q}}\{\mathbf{c}^\top \widetilde{\mathbf{y}} \} = \max_{\overline{\mathbf{c}} \in \widetilde{\mathcal{S}}} \mathbb{E}_{\mathbb{Q}}\{\overline{\mathbf{c}}^\top \widetilde{\mathbf{y}} \},
\end{equation}
where $\widetilde{\mathcal{S}}$ is a polyhedral uncertainty set corresponding to the ambiguity set $\widetilde{\mathcal{Q}} \in \{ \mathcal{Q}_1, \ldots, \mathcal{Q}_{2^{|\mathcal{L}|}}\}$ identified by the user while traversing through the network; recall Theorem~\ref{theorem 1}.

We observe that the sum $\rho_1 + \rho_2$ characterizes the total value (in percentages) by which the difference $z^*_{static} - z_{lower}$ can be reduced by leveraging our two-step validation approach. As a potential drawback, it may be the case that $\widetilde{z}_{dynamic} < z_{lower}$ and, thus, the actual upper bound for $\rho_2$ is not well-defined.  

\subsection{Test instances} \label{subsec: test instances}
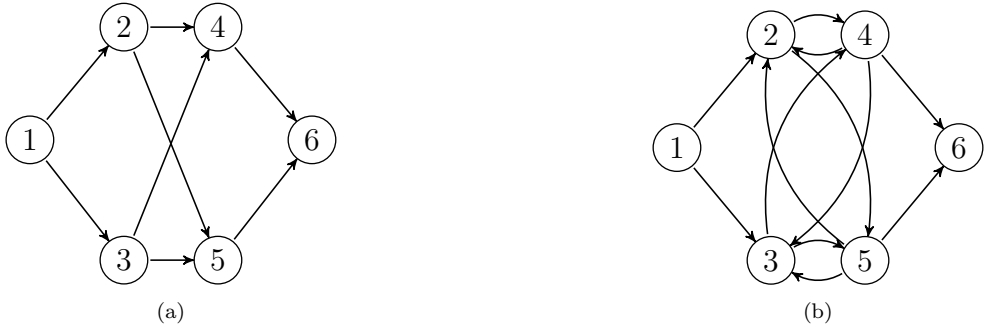
\begin{figure}
	\begin{subfigure}[b]{0.5\textwidth}
		\centering
		\begin{tikzpicture}[scale = 0.5]
		\Vertex[x=0,y=0.4]{1}
		\Vertex[x=2.5,y=3.4]{2}
		\Vertex[x=2.5,y=-2.8]{3}
		\Vertex[x=5,y=3.4]{4}
		\Vertex[x=5,y=-2.8]{5}
		\Vertex[x=7.5,y=0.4]{6}
		\tikzstyle{LabelStyle}=[fill=white,sloped]
		\tikzstyle{EdgeStyle}=[post]
		\Edge(1)(2)
		\Edge(1)(3)
		\Edge(2)(4)
		\Edge(2)(5)
		\Edge(3)(4)
		\Edge(3)(5)
		\Edge(4)(6)
		\Edge(5)(6)
		\end{tikzpicture}
		\caption{ }
		\label{fig: acyclic layered graph}
	\end{subfigure}
	\begin{subfigure}[b]{0.5\textwidth}
		\centering
		\begin{tikzpicture}[scale = 0.5]
		\Vertex[x=0,y=0]{1}
		\Vertex[x=2.5,y=3]{2}
		\Vertex[x=2.5,y=-3]{3}
		\Vertex[x=5,y=3]{4}
		\Vertex[x=5,y=-3]{5}
		\Vertex[x=7.5,y=0]{6}
		\tikzstyle{LabelStyle}=[fill=white,sloped]
		\tikzstyle{EdgeStyle}=[post, bend left]
		\Edge(2)(4)
		\Edge(4)(2)
		\Edge(3)(5)
		\Edge(5)(3)
		\Edge(2)(5)
		\Edge(5)(2)
		\Edge(3)(4)
		\Edge(4)(3)
		\tikzstyle{EdgeStyle}=[post]
		\Edge(1)(2)
		\Edge(1)(3)
		\Edge(4)(6)
		\Edge(5)(6)
		\end{tikzpicture}
		\caption{ }
		\label{fig: general layered graph}
	\end{subfigure}	
	\caption{\footnotesize \centering Acyclic (\ref{fig: acyclic layered graph}) and general (\ref{fig: general layered graph}) layered graphs with $h = 2$ intermediate layers and $r = 2$ nodes at each layer. The source and the destination nodes are given by $s = 1$ and $f = 6$, respectively.}
	\label{fig: acyclic and general graphs}
\end{figure}

\textbf{Classes of graphs.}
We consider two classes of fully connected layered graphs that are either acyclic or contain directed cycles. An acyclic layered graph is assumed to contain $h \in \mathbb{Z}_+$ intermediate layers and $r \in \mathbb{Z}_+$ nodes at each layer. The first and the last layer consist of unique nodes, which are the source and the destination nodes, respectively. General graphs are, in turn, obtained from acyclic graphs assuming that the arcs not adjacent to $s$ and $f$ can be traversed in both directions; see, e.g., Figure \ref{fig: acyclic and general graphs}. 



\looseness-1 \textbf{Nominal distribution and initial ambiguity set.} 
With respect to the nominal marginal distributions $\mathbb{Q}_a^*$, $a \in A$, we assume that the arc costs $c_a$ are governed by a standard beta distribution with parameters $\alpha_a, \beta_a \in \mathbb{R}_+$ and a support given by $[0, 1]$. 
 The parameters $\alpha_a$ and $\beta_a$ for each $a \in A$ can be defined using the mean, $m_a$, and the standard deviation, $\sigma_a$, of $c_a$, i.e.,
\begin{equation} \label{eq: parameters of beta distribution}
\begin{gathered}
\alpha_a = \frac{m_a^2(1 - m_a)}{\sigma^2_a} - m_a, \quad
\beta_a = \alpha_a(\frac{1}{m_a} - 1),
\end{gathered}
\end{equation}
see, e.g., \cite{Gupta2004}.
In particular, for each $a = (i, j) \in A$, $i < j$, we set $\sigma_{a} = 0.125$, and select $m_{a}$ uniformly at random from the interval
$$\Big(\frac{1}{2}(1 - \sqrt{1 - 4\sigma_{a}^2}), \frac{1}{2}(1 + \sqrt{1 - 4\sigma_{a}^2})\Big).$$
The latter condition guarantees that a beta distribution defined by (\ref{eq: parameters of beta distribution}) exists, i.e., $\alpha_{a}, \beta_{a} > 0$. 
 In addition, we assume that the nominal marginal distributions $\mathbb{Q}^*_{(i, j)}$ and $\mathbb{Q}^*_{(j, i)}$ are \textit{the same} for any $(i, j) \in A$. Finally, the joint distribution $\mathbb{Q}^*$ is defined as a product of marginal distributions $\mathbb{Q}^*_a$, $a \in A$. 

Next, the initial ambiguity set $\mathcal{Q}$ is constructed as follows:
\begin{equation} \label{initial ambiguity set}
\begin{gathered}
\mathcal{Q} = \Big\{\mathbb{Q} \in \mathcal{Q}_0(\mathbb{R}^{|A|}) :\; \mathbb{E}_{\mathbb{Q}} \Big\{\sum_{a \in FS_i \cup RS_i} c_a - \sum_{\{(i,j):\; (i,j) \in FS_i, \mbox { } (j,i) \in RS_i\}} c_{(i,j)} \Big\} \leq \Gamma_i
\quad \forall i \in N, \\ 
\mathbb{E}_{\mathbb{Q}}\{c_{(i, j)} - c_{(j, i)}\} = 0 \quad \forall (i, j) \in A, \;
\mathbb{Q}_a\{c_a \in [0, 1]\} = 1 \quad \forall a \in A 
\Big\}. \qquad \;
\end{gathered}
\end{equation}
\looseness-1 Specifically, the linear expectation constraints with respect to $a \in FS_i \cup RS_i$ can be thought as some budget constraints from the attacker's perspective (if the graph is general, then we account each arc $a \in FS_i \cup RS_i$ only once). 
The second linear expectation constraints in (\ref{initial ambiguity set}) indicate that the expected costs of $(i, j) \in A$ and $(j, i) \in A$ are the same; these constraints stem from the definition of $\mathbb{Q}^*$ and are assumed to be satisfied by construction. The last constraints in (\ref{initial ambiguity set}) are the support constraints.


For every $i \in N$ we calculate the parameters $\Gamma_i$ in (\ref{initial ambiguity set}) by leveraging a data set
\begin{equation} \nonumber
\widetilde{\mathbf{C}} = \Big\{(\tilde{c}_1^{(k)}, \ldots, \tilde{c}_{|A|}^{(k)})^\top, \; k \in \{1, \ldots, \widetilde{n} \} \Big\}
\end{equation} 
obtained from the nominal distribution $\mathbb{Q}^*$ and Hoeffding inequality (\ref{eq: hoeffding 2}) with some prescribed confidence level $\eta_i \in (0, 1)$. 
 Specifically, in view of Bonferroni's inequality \cite{Birge2011}, we set $\eta_i = 1 - \frac{1 - \eta}{|N|}$, $i \in N$, where $\eta \in (0, 1)$ is a required confidence level for the ambiguity set $\mathcal{Q}$. 
For convenience, we report the notations used in our test instances in Table~\ref{tab: parameters}. 

\begin{table}
	\centering
	\footnotesize
	\onehalfspacing
	
	\begin{tabular}{c| c}
		\hline 
		Parameter & Definition \\
		\hline
		$h \in \mathbb{Z}_{>0}$ & a number of intermediate layers in the graph \\
		$r \in \mathbb{Z}_{>0}$ & a number of nodes at each layer \\
		$\kappa \in (0, 1)$ & a probability that a sensor is placed at each node \\
		$\widetilde{n} \in \mathbb{Z}_{>0}$ & a number of samples in the data set $\widetilde{\mathbf{C}}$ \\
		$\widehat{n} \in \mathbb{Z}_{>0}$ & a number of samples in the data set $\widehat{\mathbf{C}}$ \\
		$\eta_i = 1 - \frac{1 - \eta}{|N|}$, $i \in N$ & confidence levels for the linear expectation constraints \\
		$\eta \in (0, 1)$ & a confidence level for the initial ambiguity set (\ref{initial ambiguity set}) \\
		$\gamma \in (0, 1)$ & a confidence level for each auxiliary constraint \\
		\hline
	\end{tabular}
	\caption{\footnotesize Summary of the notation used in the construction of test instances.}
	\label{tab: parameters}
\end{table}

Summarizing the discussion above, our construction of $\mathcal{Q}$ is somewhat stylized and not without limitations. However, in Section \ref{subsec: result for different types of networks} it is verified numerically that $z^*_{static} \neq z_{lower}$ in all generated test instances (otherwise, there is no need to apply our multi-stage problem formulation). Furthermore, we do not consider some more involved initial ambiguity sets as our focus is on the role of the auxiliary distributional constraints~(\ref{eq: list of additional constraints}).

\textbf{Auxiliary constraints.}
For simplicity, we suppose that the auxiliary constraints are restricted by the linear expectation constraints of the form (\ref{cons: auxiliary expectation}); some additional results for the auxiliary constraints (\ref{cons: auxiliary probability}) are relegated to \ref{sec: app}.
We recall that by Theorem~\ref{theorem 1} the probability constraints (\ref{cons: auxiliary probability}) can also be seen as some interval constraints with respect to expected costs. However, the difference between (\ref{cons: auxiliary probability}) and (\ref{cons: auxiliary expectation}) lies in the form of resulting expectation constraints and the data used to verify these constraints; we refer to our discussion in Section~\ref{subsec: constraints form data} and~\ref{sec: app}. 

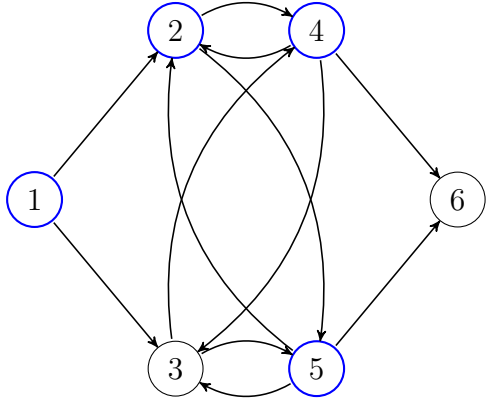
\begin{figure}
	\centering
	\begin{tikzpicture}[scale = 0.75]
	\tikzstyle{VertexStyle}=[draw = blue, circle, thick]
	\Vertex[x=2.5,y=3]{2}
	\Vertex[x=5,y=3]{4}
	\Vertex[x=5,y=-3]{5}
	\Vertex[x=0,y=0]{1}
	\tikzstyle{VertexStyle}=[draw = black, circle, thin]
	\Vertex[x=2.5,y=-3]{3}
	\Vertex[x=7.5,y=0]{6}
	\tikzstyle{LabelStyle}=[fill=white,sloped]
	\tikzstyle{EdgeStyle}=[post, bend left]
	\Edge(2)(4)
	\Edge(4)(2)
	\Edge(3)(5)
	\Edge(5)(3)
	\Edge(2)(5)
	\Edge(5)(2)
	\Edge(3)(4)
	\Edge(4)(3)
	\tikzstyle{EdgeStyle}=[post]
	\Edge(1)(2)
	\Edge(1)(3)
	\Edge(4)(6)
	\Edge(5)(6)
	\end{tikzpicture}	
	\caption{\footnotesize An example of random sensors allocation for a layered graph with $h = r = 2$. The sensors are highlighted in blue.}
	\label{fig: general layered graph with sensors}
\end{figure}

We assume that each node $i \in N$ is equipped with a sensor with some fixed probability $\kappa \in (0, 1)$; see, e.g., Figure \ref{fig: general layered graph with sensors}. 
 Furthermore, following Example \ref{example 2}, we suggest three types of auxiliary expectation constraints associated with arcs in the forward direction from $s$ to $f$:
\begin{itemize}
	\item \textit{individual constraints} corresponding to the adjacent sensors, e.g., constraints for the expected cost of $(1,2)$, $(2,4)$ and $(2,5)$ in Figure \ref{fig: general layered graph with sensors}. 
	\item \textit{difference constraints}, e.g., constraints with respect to the difference of expected costs of $(3, 4)$ and $(3,5)$ in Figure \ref{fig: general layered graph with sensors}. 
	\item \textit{sum constraints}, e.g., constraints with respect to the sum of expected costs of $(3, 4)$ and $(3,5)$ in Figure \ref{fig: general layered graph with sensors} (we use these constraints, if at least one of the arcs can be traversed in both directions). 
\end{itemize}

Formally, a constraint of the form (\ref{cons: auxiliary expectation}) related to some node $i \in N \setminus \{f\}$ can be defined as:
\begin{equation} \label{eq: sensor constraints}
\mathbb{E}_{\mathbb{Q}} \Big \{\sum_{a \in FS_i} p_a c_a \Big \} = \sum_{a \in FS_i} p_a \overline{c}_a \leq p_0,
\end{equation}
where $\overline{c}_a = \mathbb{E}_{\mathbb{Q}_a}\{c_a\}$, the coefficients $p_a \in \{-1, 0, 1\}$ depend on the type of this constraint and the threshold $p_0 \in \mathbb{R}$ is, e.g., defined as some average value of the sum $\sum_{a \in FS_i} p_a \overline{c}_a$ over the initial polyhedron $\mathcal{S}_0$: 
\begin{equation} \nonumber
p_0 := \frac{1}{2} \Big(\min_{\overline{\mathbf{c}} \in \mathcal{S}_0} \sum_{a \in FS_i} p_a \overline{c}_a + \max_{\overline{\mathbf{c}} \in \mathcal{S}_0} \sum_{a \in FS_i} p_a \overline{c}_a\Big). 
\end{equation}
We collect a required number of the auxiliary constraints by selecting a random node and a random constraint associated with this node. If the above procedure does not allow to construct the required number of constraints, then the sensors allocation is updated. 

For convenience, we generate a new data set 
\begin{equation} \nonumber
\widehat{\mathbf{C}} := \Big\{(\hat{c}_1^{(k)}, \ldots, \hat{c}_{|A|}^{(k)})^\top, \; k \in \{1, \ldots, \widehat{n} \} \Big\}
\end{equation}
\looseness-1 for the constraint verification procedure; recall that we used another data set $\widetilde{\mathbf{C}}$ to construct the initial ambiguity set $\mathcal{Q}$. This splitting of data is rather technical and reflects the fact that the data for the constraint verification procedure can be collected in an online manner and, thus, it is not necessarily available at the beginning of the game. 

The auxiliary constraints are verified via Hoeffding inequality~(\ref{eq: hoeffding 2}) with some fixed confidence level $\gamma \in (0, 1)$ (which is assumed to be the same for all auxiliary constraints). As outlined before, in order to verify a constraint of the form (\ref{eq: sensor constraints}) related to some node~$i~\in~N~\setminus~\{f\}$, we~use only random samples (or their linear combinations) associated with the arcs $a \in FS_i$. If the user fails to verify some auxiliary constraint using $\widehat{n}$ samples, then this constraint is verified by force as follows:
\begin{itemize}
	\item if the constraint is an individual or a sum constraint, then it is assumed to be violated (it models the worst-case scenario for the user);
	\item if the constraint is a difference constraint, then it is assumed to be satisfied with probability 0.5.
\end{itemize}
Summarizing the discussion above, we expect that for some fixed confidence level $\gamma$ the quality of the constraint verification procedure, (\ref{eq: relative gap 2}), increases with the increase of $\widehat{n}$.

\begin{table}
	\footnotesize
	\centering
	\onehalfspacing
	
	\begin{tabular}{c | c c c }
		\multirow{2}{*}{Number of constraints} & \multicolumn{3}{c}{Acyclic graphs} \\\cline{2-4}
		& $\rho_1$ (MAD) in $\%$ & $\rho_2$ (MAD) in $\%$ & average time (MAD) in sec. \\\hline
		$|\mathcal{L}| = 1$ & 2.0 (3.5) 
		& 1.5 (2.8) & 0.04 (0.02) \\
		$|\mathcal{L}| = 2$ & 4.8 (7.2) & 5.8 (9.1) & 0.07 (0.02) \\		
		$|\mathcal{L}| = 3$ & 6.3 (8.9) & 8.3 (11.7) & 0.17 (0.07) \\
		$|\mathcal{L}| = 4$ & 8.3 (10.0) & 11.3 (14.1) & 0.67 (0.33) \\
		$|\mathcal{L}| = 5$ & 10.2 (11.6) & 12.2 (15.2) & 3.8 (2.85) \\
		\hline
	\end{tabular}
	\caption{\footnotesize Let $h = r = 3$, $\widehat{n} = \widetilde{n} = 60$ and assume that the graph is acyclic. We report the average relative gaps (\ref{eq: relative gap 1}) and (\ref{eq: relative gap 2}) and the average running times in seconds with mean absolute deviations (in brackets) over $50$ random test instances.}
	\label{tab: results for acyclic}
\end{table}

\begin{table}
	\footnotesize
	\centering
	\onehalfspacing
	
	\begin{tabular}{c | c c c}
		\multirow{2}{*}{Number of constraints} & \multicolumn{3}{c}{General graphs} \\\cline{2-4}
		& $\rho_1$ (MAD) in $\%$ & $\rho_2$ (MAD) in $\%$ & average time (MAD) in sec. \\\hline
		$|\mathcal{L}| = 1$ & 1.8 (3.2) & 3.3 (5.9) & 0.08 (0.05) \\ 			
		$|\mathcal{L}| = 2$ & 3.6 (6.0) & 4.4 (7.3) & 0.28 (0.10) \\
		$|\mathcal{L}| = 3$ & 4.9 (7.1) & 9.1 (13.1) & 10.1 (6.71) \\
		$|\mathcal{L}| = 4$ & - & - & $>600$ \\
		$|\mathcal{L}| = 5$ & - & - & $>600$ \\
		\hline
	\end{tabular}
	\caption{\footnotesize Let $h = r = 3$, $\widehat{n} = \widetilde{n} = 60$ and assume that the graph is general. We report the average relative gaps (\ref{eq: relative gap 1}) and (\ref{eq: relative gap 2}) and the average running times in seconds with mean absolute deviations (in brackets) for $50$ random test instances.}
	\label{tab: results for general}
\end{table}

\subsection{Results for different network types} \label{subsec: result for different types of networks}
In the numerical experiments we explore how the quality of adaptive decisions scales in the number of auxiliary constraints, the size of the network and the sample sizes, $\widehat{n}$ and $\widetilde{n}$; recall Table \ref{tab: parameters}. In particular, for different parameter settings we compute the average relative gaps (\ref{eq: relative gap 1}) and (\ref{eq: relative gap 2}) and the average running time with mean absolute deviations for $50$ test instances. Also, in all experiments we set $\eta = 0.95$, $\gamma = 0.95$ and $\kappa = 0.5$. 

\textbf{Dependence on the number of auxiliary constraints.} Let $h = r = 3$ and $\widehat{n} = \widetilde{n} = 60$. For each test instance we increase the number of auxiliary constraints, $|\mathcal{L}|$, from $1$ to $5$. The results  for acyclic and general graphs are reported, respectively, in Tables \ref{tab: results for acyclic} and~\ref{tab: results for general}.

 The observations from our computational results can be summarized as follows: 
\begin{itemize}
	\item Augmenting the set of auxiliary
	constraints provides more information for the user and, thus, improves the quality of both adaptive decisions and the constraint verification procedure. For example, for $|\mathcal{L}| = 5$ and acyclic graphs the total profit $\rho_1 + \rho_2$ obtained by the user is about $22\%$ on average.
	\item The multi-stage problem (\ref{multi-stage DRSPP}) becomes more computationally complex for general graphs and with the increase of $|\mathcal{L}|$. This fact is rather intuitive as the number of variables and constraints in (\ref{multi-stage DRSPP}) is exponential in $|\mathcal{L}|$ and increases for general graphs. 
\end{itemize}

\begin{table}
	\footnotesize
	\centering
	\onehalfspacing
	
	\begin{tabular}{c| c| c c c }
		\multirow{2}{*}{$\#$ samples in $\widetilde{\mathbf{C}}$} & \multirow{2}{*}{$\#$ samples in $\widehat{\mathbf{C}}$} & \multicolumn{3}{c}{Acyclic graphs} \\\cline{3-5} 
		& & $\rho_1$ (MAD) in $\%$ & $\rho_2$ (MAD) in $\%$ & average time (MAD) in sec \\\hline
		$\widetilde{n} = 20$ & \multirow{5}{*}{$\widehat{n} = 60$} & 2.0 (3.6) & 22.7 (29.1) & 2.49 (1.80) \\ 			
		$\widetilde{n} = 40$ & & 7.9 (9.7) & 16.3 (18.5) & 3.43 (2.46) \\
		$\widetilde{n} = 60$ & & 10.2 (11.6) & 12.2 (15.2) & 3.8 (2.85) \\
		$\widetilde{n} = 80$ & & 11.6 (11.8) & 9.5 (11.9) & 5.45 (5.11) \\
		$\widetilde{n} = 100$ & & 12.2 (11.6) & 12.2 (13.2) & 3.89 (2.95) \\	
		\hline
		\multirow{5}{*}{$\widetilde{n} = 60$} & $\widehat{n} = 20$ & \multirow{5}{*}{10.2 (11.6)} & 9.9 (12.8) & \multirow{5}{*}{3.8 (2.85)}\\ 			
		& $\widehat{n} = 40$ & & 11.4 (14.9) & \\
		& $\widehat{n} = 60$ & & 12.2 (15.2) & \\
		& $\widehat{n} = 80$ & & 13.1 (15.6) & \\
		& $\widehat{n} = 100$ & & 13.2 (16.3) & \\	
		\hline
		
	\end{tabular}
	\caption{\footnotesize Let $h = r = 3$, $|\mathcal{L}| = 5$ and assume that the graph is acyclic. We report the average relative gaps (\ref{eq: relative gap 1}) and (\ref{eq: relative gap 2}) and the average running times in seconds with mean absolute deviations (in brackets) for $50$ random test instances.}
	\label{tab: results for n0 acyclic}
\end{table}


 
 \textbf{Dependence on the sample size.} In this experiment we consider acyclic graphs and explore the value of parameters $\widehat{n}$ and $\widetilde{n}$. For each test instance and different values of $\widehat{n}$ and $\widetilde{n}$ we preserve the same nominal distribution $\mathbb{Q}^*$ and fix a particular placement of sensors. The results for $h = r = 3$ and $|\mathcal{L}| = 5$ are reported in Table \ref{tab: results for n0 acyclic}. 

We make the following observations:
\begin{itemize}
	\item (\ref{eq: relative gap 1}) and (\ref{eq: relative gap 2}) tend to increase and decrease, respectively, as functions of~$\widetilde{n}$. It can be argued that for sufficiently small $\widetilde{n}$ we have $z^*_{static} \approx z_{lower}$. Our hypothesis is that in this case dynamic decisions cannot significantly improve the value of $z^*_{static}$, but the constraint verification procedure allows to reimburse these losses. 
	\item (\ref{eq: relative gap 2}) increases as a function of $\widehat{n}$; see the second part of Table \ref{tab: results for n0 acyclic}. This fact is quite intuitive, since with the increase of $\widehat{n}$ the user may verify more constraints correctly.
	\item The solution times for the MIP formulation (\ref{MIP reformulation}) do not depend on $\widehat{n}$ by construction and slightly depend on $\widetilde{n}$; see the last column of Table \ref{tab: results for n0 acyclic}.
\end{itemize}

\textbf{Dependence on the size of the network.} In our last numerical experiment we examine the value of the network's size. 
Let $\widehat{n} = \widetilde{n} = 60$ and assume that $|\mathcal{L}| = 5$ and $|\mathcal{L}| = 3$ for acyclic and general graphs, respectively. For each class of graphs we test different values of $h \in \{2, 3, 4\}$ and $r \in \{2, 3, 4\}$ and report the results in Tables \ref{tab: results for size acyclic} and \ref{tab: results for size general}.

We make the following observations:
\begin{itemize}
	\item The solution times increase in both $h$ and $r$, while the average relative gaps (\ref{eq: relative gap 1}) and~(\ref{eq: relative gap 2}) tend to decrease. Indeed, the size of the MIP reformulation (\ref{MIP reformulation}) increases as a function of $h$ and $r$, while the role of each particular constraint in $\mathcal{L}$ becomes less pronounced for networks of a larger size.
	\item The average relative gap (\ref{eq: relative gap 1}) tends to zero faster with the increase of $r$ rather than $h$. In this regard, we note that the larger $r$, the less constraints are potentially met by the user while traversing through the network.
\end{itemize}


\textbf{Summary.} 
It turns out that the user may often gain some profit by leveraging dynamic decisions instead of static ones. However, the magnitude of this profit depends on the size of the graph, the form of initial ambiguity set as well as the parameters of the constraint verification procedure. Admittedly, the primary application of our approach is to the networks of a relatively small size, in which solutions with a sufficiently high quality can be derived within a reasonable time. It can be also argued that most of existing solution approaches to multi-stage problems with binary recourse decisions can only be applied to instances of a moderate size; see, e.g., the studies in \cite{Bertsimas2016, Bertsimas2015, Postek2016} and the references therein.

\begin{table}
	\footnotesize
	\centering
	\onehalfspacing
	
	\begin{tabular}{c c | c c c}
		\multirow{2}{*}{$\#$ layers} & \multirow{2}{*}{$\#$ nodes at each layer} & \multicolumn{3}{c}{Acyclic graphs} \\\cline{3-5} 
		& & $\rho_1$ (MAD) in $\%$ & $\rho_2$ (MAD) in $\%$ & average time (MAD) in sec \\\hline				
		$h = 2$ & \multirow{3}{*}{$r = 3$} & 22.0 (19.8) & 19.5 (24.3) & 0.42 (0.21) \\
		$h = 3$ & & 10.2 (11.6) & 12.2 (15.2) & 3.8 (2.85)\\
		$h = 4$ & & 7.6 (9.3) & 8.6 (11.5) & 34.46 (30.34) \\
		\hline
		\multirow{3}{*}{$h = 3$} & $r = 2$ & 19.8 (16.3) & 34.1 (26.2) & 0.19 (0.05) \\
		& $r = 3$ & 10.2 (11.6) & 12.2 (15.2) & 3.8 (2.85) \\
		& $r = 4$ & 3.3 (5.7) & 11.6 (15.9) & 44.13 (47,27) \\
		\hline
	\end{tabular}
	\caption{\footnotesize Let $|\mathcal{L}| = 5$, $\widehat{n} = \widetilde{n} = 60$ and assume that the graph is acyclic. We report the average relative gaps (\ref{eq: relative gap 1}) and (\ref{eq: relative gap 2}) and the average running times in seconds with mean absolute deviations (in brackets) for $50$ random test instances.}
	\label{tab: results for size acyclic}
\end{table}

\begin{table}
	\footnotesize
	\centering
	\onehalfspacing
	
	\begin{tabular}{c c | c c c}
		\multirow{2}{*}{$\#$ layers} & \multirow{2}{*}{$\#$ nodes at each layer} & \multicolumn{3}{c}{General graphs} \\\cline{3-5} 
		& & $\rho_1$ (MAD) in $\%$ & $\rho_2$ (MAD) in $\%$ & average time (MAD) in sec \\\hline				
		$h = 2$ & \multirow{3}{*}{$r = 3$} & 16.6 (18.2) & 10.4 (15.6) & 0.72 (0.27) \\
		$h = 3$ & & 4.9 (7.1) & 9.1 (13.1) & 10.1 (6.71) \\
		$h = 4$ & & 2.5 (4.2) & 0.9 (1.7) & 99.62 (88.53) \\
		\hline
		\multirow{3}{*}{$h = 3$} & $r = 2$ & 16.0 (16.9) & 9.9 (13.7) & 0.47 (0.12)\\
		& $r = 3$ & 4.9 (7.1) & 9.1 (13.1) & 10.1 (6.71) \\
		& $r = 4$ & 1.1 (2.1) & 1.7 (3.1) & 64.54 (43.12) \\
		\hline
	\end{tabular}
	\caption{\footnotesize Let $|\mathcal{L}| = 3$, $\widehat{n} = \widetilde{n} = 60$ and assume that the graph is general. We report the average relative gaps (\ref{eq: relative gap 1}) and (\ref{eq: relative gap 2}) and the average running times in seconds with mean absolute deviations (in brackets) for $50$ random test instances.}	\label{tab: results for size general}
\end{table}

\section{Conclusion} \label{sec: conclusion}

\looseness-1In this paper we consider a dynamic version of the shortest path problem, where the cost vector is subject to distributional uncertainty. We formulate the problem in terms of an ambiguity-averse multi-stage network game between a user and an attacker. The user aims at minimizing its cumulative expected loss by traversing between two fixed nodes in the network, while the attacker maximizes the user's objective function by selecting a distribution of arc costs from an ambiguity set of candidate distributions. In contrast to the one-stage formulation, both the user and the attacker are able to adjust their decisions at particular nodes of the user's path. 

Following the related one-stage formulation, we suppose that the family of probability distributions is formed by individual probability constraints for particular arcs and linear expectation constraints for prescribed subsets of arcs. Also, while traversing through the network, the user may verify some auxiliary distributional constraints associated with the arcs emanated from the current user's position. Specifically, we assume that the user forms a list of auxiliary distributional constraints at the beginning of the game, whereas the attacker forms a list of responses that are consistently revealed to the user. 

We design two classes of non-anticipativity constraints (for acyclic and general graphs, repsectively) enforcing that the user's decision at some particular node cannot depend on future attacker's responses. Furthermore, by using some properties of the related one-stage formulation and linear programming duality the multi-stage problem is reformulated as a one potentially large linear mixed-integer programming problem.

From the application perspective, we illustrate that the auxiliary distributional constraints can be constructed and verified by using some information from Bluetooth sensors. Finally, we conduct a numerical study where the one- and multi-stage problem formulations are compared with respect to several classes of synthetic network instances. It turns out that using adaptive decisions is practically relevant only for networks of a relatively small size, in which high quality adaptive decisions can be obtained within a reasonable time. 

With respect to future research directions, it would be interesting to explore some other types of objective criteria and ambiguity sets; see, e.g., \cite{Hanasusanto2016}, in the context of our multi-stage shortest path problem. At the same time, it can be argued that another types of ambiguity sets imply the need for more advanced solution techniques as our partitions of the initial ambiguity set are based on some linearity properties of the auxiliary distributional constraints. As another possible option one may consider randomized decisions from the user's perspective; see, e.g., the related studies in \cite{Delage2019, Delage2022}.

\textbf{Acknowledgments.} The author would like to thank Dr. Oleg Prokopyev for his helpful comments and suggestions. The article was prepared within the framework of the Basic Research Program at the National Research University Higher School of Economics (Sections~ 1-2). The research is funded by RSF grant №22-11-00073~(Sections 3-5). The author is thankful to the associate editor and two anonymous referees for their constructive comments that allowed to greatly improve the paper. 

\textbf{Conflict of interest:} None.

\textbf{Data availability statement:} The author confirms that all data generated or analysed during this study are included in this published article.
\bibliographystyle{apa}
\bibliography{bibliography}

\appendix 

\section{Supplementary material} \label{sec: app}
In this appendix we consider the ambiguity set (\ref{initial ambiguity set}) with the auxiliary probability constraints of the form:
\begin{equation} \label{eq: axuxiliary probability constraint experiments}
\mathbb{Q}_a \{c_a \in [0.5, 1]\} \leq \widetilde{q}_a
\end{equation}
with $\widetilde{q}_a = \widetilde{q}$ for each $a \in A$; recall Section \ref{subsec: constraints form data}. As before, we assume that each node $i \in N$ is equipped with a sensor with probability $\kappa = 0.5$ and the constraints (\ref{eq: axuxiliary probability constraint experiments}) correspond to the arcs $a \in A$ with the adjacent sensors. 
In all experiments we set $h = r = 3$, $\widehat{n} = \widetilde{n} = 60$, $\alpha = \gamma = 0.95$ and assume that the graph is acyclic; recall Table \ref{tab: parameters}. 

For the first experiment we set $\widetilde{q} = 0.5$ and consider the relative gaps (\ref{eq: relative gap 1}) and (\ref{eq: relative gap 2}), as well as the running time for solving the MIP reformulation (\ref{MIP reformulation}) as a function of $|\mathcal{L}|$. The results are reported in Table \ref{tab: results for acyclic 2}. For the second experiment we set $|\mathcal{L}| = 5$ and explore how the quality of adaptive decisions scales in the parameter $\widetilde{q}$ of the auxiliary constraint (\ref{eq: axuxiliary probability constraint experiments}). For each test instance and different values of $\widetilde{q}$ we preserve the same data set (\ref{eq: data set}) and the same initial ambiguity set (\ref{initial ambiguity set}). The results are reported in Table \ref{tab: results dependence on q}.

We make the following observations:
\begin{itemize}
	\item Comparing the obtained results in Tables \ref{tab: results for acyclic} and \ref{tab: results for acyclic 2} we observe that the quality of adaptive solutions decreases whenever the probability constraints (\ref{cons: auxiliary probability}) are used instead of the expectation constraints (\ref{cons: auxiliary expectation}). In fact, if the constraint (\ref{eq: axuxiliary probability constraint experiments}) is satisfied, then 
	\begin{equation} \label{eq: app expectation 1}
	\mathbb{E}_{\mathbb{Q}_a}\{c_a\} \leq 0.5 \times (1 - \widetilde{q}) + 1 \times \widetilde{q} = 0.5 \; \widetilde{q} + 0.5.
	\end{equation}
	Otherwise, we have 
	\begin{equation} \label{eq: app expectation 2}
	\mathbb{E}_{\mathbb{Q}_a}\{c_a\} \geq 0 \times (1 - \widetilde{q}) + 0.5 \times \widetilde{q} = 0.5 \;\widetilde{q}.
	\end{equation}
	This observation implies that the user receives less information in terms of expected costs by using (\ref{eq: axuxiliary probability constraint experiments}) instead of an individual linear expectation constraint. 
	\item Both relative gaps (\ref{eq: relative gap 1}) and (\ref{eq: relative gap 2}) tend to decrease whenever $\widetilde{q}$ approaches to its lower or upper bounds. Indeed, if $\widetilde{q}$ close to its bounds, then verification of some constraints of the form (\ref{eq: axuxiliary probability constraint experiments}) may not improve the user's objective function value.
	\item The constraints (\ref{eq: axuxiliary probability constraint experiments}) may handle interval-censored observations, which is not the case for the auxiliary expectation constraints (\ref{eq: sensor constraints}).
\end{itemize}

\begin{table}
	\footnotesize
	\centering
	\onehalfspacing
	\begin{tabular}{c | c c c }
	\multirow{2}{*}{Number of constraints} & \multicolumn{3}{c}{Acyclic graphs} \\\cline{2-4}
	& $\rho_1$ (MAD) in $\%$ & $\rho_2$ (MAD) in $\%$ & average time (MAD) in sec. \\\hline
	$|\mathcal{L}| = 1$ & 1.1 (2.1) 
	& 1.8 (3.1) & 0.03 (0.01) 
	 
	 \\
	$|\mathcal{L}| = 2$ & 1.7 (2.9) & 6.3 (8.8) & 0.06 (0.02) \\		
	$|\mathcal{L}| = 3$ & 2.0 (3.4) & 10.2 (10.7) & 0.18 (0.07) \\
	$|\mathcal{L}| = 4$ & 2.4 (3.9) & 12.0 (11.7) & 0.69 (0.38) \\
	$|\mathcal{L}| = 5$ & 2.9 (4.4) & 13.6 (12.9) & 4.03 (3.05) \\
	\hline
\end{tabular}
\caption{\footnotesize Let $h = r = 3$, $\widehat{n} = \widetilde{n} = 60$, $\overline{q} = 0.5$ and assume that the graph is acyclic. We report the average relative gaps (\ref{eq: relative gap 1}) and (\ref{eq: relative gap 2}) and the average running times in seconds with mean absolute deviations (in brackets) for $50$ random test instances.}
\label{tab: results for acyclic 2}
\end{table}

\begin{table}
	\footnotesize
	\centering
	\onehalfspacing
	\begin{tabular}{c | c c c }
		\multirow{2}{*}{Parameter $\widetilde{q}$} & \multicolumn{3}{c}{Acyclic graphs} \\\cline{2-4}
		& $\rho_1$ (MAD) in $\%$ & $\rho_2$ (MAD) in $\%$ & average time (MAD) in sec. \\\hline
		$\widetilde{q} = 0.1$ & 0.5 (0.9) 
		& 0.2 (0.3) & 3.06 (1.98) 
		\\
		$\widetilde{q} = 0.3$ & 1.6 (2.6) & 15.8 (16.1) & 4.62 (3.89) \\		
		$\widetilde{q} = 0.5$ & 2.9 (4.4) & 13.6 (12.9) & 4.03 (3.05) \\
		$\widetilde{q} = 0.7$ & 3.4 (4.4) & 8.3 (8.0) & 4.78 (4.16) \\
		$\widetilde{q} = 0.9$ & 1.5 (1.8) & 2.7 (3.4) & 4.84 (3.85) \\
		\hline
	\end{tabular}
	\caption{\footnotesize Let $h = r = 3$, $\widehat{n} = \widetilde{n} = 60$, $|\mathcal{L}| = 5$ and assume that the graph is acyclic. We report the average relative gaps (\ref{eq: relative gap 1}) and (\ref{eq: relative gap 2}) and the average running times in seconds with mean absolute deviations (in brackets) for $50$ random test instances.}
	\label{tab: results dependence on q}
\end{table}


In conclusion, we point out that the results of this section are rather consistent with the results of Section \ref{subsec: result for different types of networks}. That is, the quality of adaptive decisions and the computational complexity of the MIP formulation (\ref{MIP reformulation}) increase in the number of auxiliary distributional constraints. 

\end{document}